\documentclass[12pt,a4paper]{article}
\usepackage[left=2.25 cm,right=2.25cm,top=3 cm,bottom=3 cm]{geometry}
\usepackage[T1]{fontenc}

\usepackage[cp1250]{inputenc}

\usepackage{graphicx}
\usepackage{amsmath,amsthm,amsfonts,amssymb}

\usepackage{graphics}
 \input xypic
 \input epsf
 \input xy
 \xyoption{all}

\def\Ker{\mathop{\hbox{Ker}}}

\newtheorem{theorem}{Theorem}[section]
\newtheorem{definition}[theorem]{Definition}

\newtheorem{corollary}[theorem]{Corollary}

\newtheorem{example}[theorem]{Example}
\newtheorem{lemma}[theorem]{Lemma}

\newtheorem{proposition}[theorem]{Proposition}
\newtheorem{remark}[theorem]{Remark}

\newcommand{\mS}{\mathbb S}

\newcommand{\mC}{\mathbb C}

\newcommand{\mN}{\mathbb N}

\newcommand{\mZ}{\mathbb Z}

\newcommand{\mR}{\mathbb R}

\newcommand{\be}{\begin{eqnarray}}
\newcommand{\ee}{\end{eqnarray}}
\newcommand{\bd}{\begin{definition}}
\newcommand{\ed}{\end{definition}}

\newcommand{\br}{\begin{remark}}
\newcommand{\er}{\end{remark}}

\newcommand{\bt}{\begin{tabular}}
\newcommand{\et}{\end{tabular}}

\def\Ker{\mathop{\hbox{Ker}}}

\def\sp{\mathop{\mathfrak{sp}}\nolimits}
\def\sl{\mathop{\mathfrak{sl}}\nolimits}
\def\mp{\mathop{\mathfrak{mp}}\nolimits}

\def\Id{\mathop{\rm Id}\nolimits}

\begin{document}
\title{Basic aspects of symplectic
   Clifford analysis for the symplectic Dirac operator}
\author{Hendrik De Bie, Marie Holíková, Petr Somberg}
\date{}
\maketitle

\abstract 
In the present article we study basic aspects 
of the symplectic version of Clifford analysis 
associated to the symplectic Dirac operator. Focusing mostly on the symplectic vector
space of real dimension $2$, this involves the analysis of first order symmetry 
operators, symplectic Clifford-Fourier transform, reproducing kernel for the 
symplectic Fischer product and the construction of bases of symplectic monogenics 
for the symplectic Dirac operator.

{\bf Key words:} Symplectic Dirac operator, Symmetry operators, Reproducing kernel, Fischer product,
Bases of symplectic monogenics.
  
{\bf MSC classification:} 53C27, 53D05, 81R25

\endabstract

\tableofcontents


\section{Introduction}

 Harmonic analysis is a fruitful concept built on the analysis of 
function spaces equipped with a Lie (finite, discrete, etc) group 
action. A key organizing principle in analyzing function spaces 
and decomposing them into simple building blocks is the notion of 
intertwining (differential, integral) operators. The basic example
related to orthogonal symmetry are quadratic spaces like euclidean
space, sphere or hyperbolic space and the Laplace operator acting on 
scalar valued smooth functions. Another developed example concerns
quadratic spaces, smooth functions valued in the spinor space or
the Clifford algebra and the orthogonally equivariant Dirac operator
(collectively known as the orthogonal Clifford analysis.)

 In the present article we focus on a similar structure: we consider 
the symplectic symmetry instead of the orthogonal group, as for 
function spaces we consider the symplectic spinors (smooth functions 
valued in the Segal-Shale-Weil representation) and the symplectic 
Dirac operator as a symplectic intertwining differential operator.
Namely, focusing mostly on the real dimension $2$, we develop the 
basics of symplectic Clifford analysis including the analysis of
first order symmetry operators, symplectic Clifford-Fourier transform,
reproducing kernel for the symplectic Fischer product and the construction
of bases for symplectic monogenics for the symplectic Dirac operator.

Let us briefly indicate the structure of our article. Section $2$ 
contains a review of the metaplectic Howe duality, \cite{bss}, a 
concise way of describing the space of symplectic spinors through
the invariant theory of the metaplectic Lie algebra $\mp(2n,\mR)$.   
Section $3$ treats the concept of symmetry differential operators
of the symplectic Dirac operator in general even dimensions, and then 
specializes to several explicit problems both in real and complex variables. In Section $4$, we turn our attention to 
the action of the generator of the Weyl group associated to 
$\mp(2,\mR)$, giving rise
to the symplectic Clifford-Fourier transform. In particular, we 
introduce the operator of metaplectic harmonic oscillator and find 
its spectral decomposition into eigenspaces.
Section $5$ is devoted to the symplectic Fischer product and the 
construction of its reproducing kernel. In Section $6$, we construct 
some explicit bases for the space of symplectic monogenics and 
prove several useful characterizing properties

To summarize, our work is clearly the first attempt aiming to uncover fundamental 
analytical properties of the symplectic Dirac operator. The generalization
of our results to a symplectic space of arbitrary dimension or a
proper formulation of an analogue of the Cauchy-Kovalevskaya theorem
are still missing cornerstones of symplectic 
Clifford analysis.



\section{Representation theory and symplectic spinors}

Let $(\mR^{2n},\omega)$ be the symplectic vector space with coordinates 
$x_1,\ldots, x_{n} $, $y_1,\ldots , y_n$, 
and coordinate vector fields 
$\partial_{x_1},\ldots,\partial_{x_n},\partial_{y_1}, \ldots ,\partial_{y_{n}}$ 
or, equally, a symplectic frame $e_{1},\ldots,e_{2n}$, fulfilling 
\begin{eqnarray}\label{omega}
\omega(e_{j},e_{n+j})=1,\quad \omega(e_{n+j},e_{j})=-1,\quad j=1,\ldots,n
\end{eqnarray}
and zero otherwise.

The symplectic Lie algebra $\sp(2n,\mR)$ has the matrix realization 
given by the span of
\begin{eqnarray*}
&& X_{jk}=E_{j,k}-E_{n+k,n+j},\\
&& Y_{jj}=E_{j,n+j},\\
&& Y_{jk}=E_{j,n+k}+E_{k,n+j} \text{ for } j\neq k,\\
&& Z_{jj}=E_{n+j,j},\\
&& Z_{jk}=E_{n+j,k}+E_{n+k,j} \text{ for } j\neq k,
\end{eqnarray*}
where $j,k= 1,\ldots, n$, and $E_ {k,j}$ is the $2n\times 2n$ 
matrix with $1$ on the intersection of the 
$k$-th row and the $j$-th column and zero otherwise. 
The representation of $\sp(2n,\mR)$ on the symmetric algebra of $\mR^{2n}$,
$S^*(\mR^{2n},\mC)=\mC[x_1,\ldots,x_{n},y_1,\ldots,y_n]$, is given by
\begin{eqnarray}
&& X_{jk}=x_j\partial_{x_k}-y_{k}\partial_{y_j}, \nonumber \\
&& Y_{jj}=x_{j}\partial_{y_j},\nonumber\\
&& Y_{jk}=x_{j}\partial_{y_k}+x_{k}\partial_{y_j} \text{ for } j\neq k,\nonumber\\
&& Z_{jj}=y_{j}\partial_{x_j},\nonumber\\
&& Z_{jk}=y_{j}\partial_{x_k}+y_{k}\partial_{x_j} \text{ for } j\neq k. \label{sarepr}
\end{eqnarray}
\begin{definition}\label{cliffalgdef}
The symplectic Clifford algebra $Cl_s(\mR^{2n},\omega)$ on $(\mR^{2n},\omega)$
with a basis $\{e_1,\ldots,e_{2n}\}$ is an associative unital algebra over $\mC$,
given by the quotient of the tensor algebra $T(e_1,\ldots,e_{2n})$ by a
two-sided ideal $ I\subset T(e_1,\ldots,e_{2n})$ generated by
$$
v\cdot w-w\cdot v=- \mathrm{i} \,\omega (v,w)
$$ 
for all $v,w\in\mR^{2n}$ and $\mathrm{i} \in\mC$ the complex unit.
Namely, the relations 
$e_j\cdot e_k -e_k \cdot e_j=-\mathrm{i}\omega (e_j,e_k)$
for the basis $\{e_1,\ldots,e_{2n}\}$ hold true.
\end{definition}
The symplectic Clifford algebra $Cl_s(\mR^{2n},\omega)$ is isomorphic 
to the Weyl algebra $W_{2n}$ of complex valued algebraic 
differential operators on $\mR^n$, and the symplectic Lie algebra 
$\sp(2n,\mR)$ can be realized as a subalgebra of $W_{2n}$. 

In particular, $W_{2n}$ is an associative algebra generated by 
$\{q_1,\ldots,q_{n},\partial_{q_1},\ldots,\partial_{q_n}\}$, the multiplication 
operator by $q_j$ and differentiation $\partial_{q_j}$, $j=1,\ldots,n$, and 
the symplectic Lie algebra $\sp(2n,\mR)\subset W_{2n}$ is 
\begin{eqnarray}
&& X_{jk}=q_k \partial_{q_j}+\frac{1}{2}\delta_{j,k},\nonumber \\
&& Y_{jj}=-\frac{\mathrm{i} }{2}\partial_{q_{j}}^2,\nonumber \\
&& Y_{jk}=\mathrm{i}  \partial_{q_j} \partial_{q_k} \text{ for } j\neq k,\nonumber \\
&& Z_{jj}=-\frac{\mathrm{i} }{2}{q_j}^2,\nonumber \\
&& Z_{jk}=\mathrm{i} {q_j}{q_k} \text{ for } j\neq k. \label{sswrepr}
\end{eqnarray}
We denote by ${\fam2 S}({\mathbb R^n})$ the space of Schwartz functions on ${\mathbb R}^n$.
The representation of $\sp(2n,\mR)$ on the space of polynomial symplectic spinors 
$\mathrm{Pol}(\mR^{2n},\mC)\otimes {\fam2 S}({\mathbb R^n})$ is given by the combination 
of the dual representation to \eqref{sarepr} and the representation \eqref{sswrepr},
\begin{eqnarray}\label{ssrepr}
&& X_{jk}= - x_j\partial_{x_k}+y_{k}\partial_{y_j}+q_k \partial_{q_j}+\frac{1}{2}\delta_{j,k}, \nonumber\\
&& Y_{jj}=- x_{j}\partial_{y_j}-\frac{\mathrm{i} }{2} \partial_{q_{j}}^2,\nonumber\\
&& Y_{jk}=x_{k}\partial_{y_j}+x_{j}\partial_{y_k}+\mathrm{i}  \partial_{q_j} \partial_{q_k} \text{ for } j\neq k,\nonumber\\
&& Z_{jj}=- y_{j}\partial_{x_j}-\frac{\mathrm{i} }{2}{q_j}^2,\nonumber\\
&& Z_{jk}=y_{k}\partial_{x_j}+y_{j}\partial_{x_k}+\mathrm{i}  {q_j} {q_k} \text{ for } j\neq k.
\end{eqnarray}
The action of the basis elements $e_1,\ldots,e_{2n}$ on a symplectic spinor 
$\varphi\in {\mathcal S}({\mathbb R}^n)$ is given by 
\begin{eqnarray}
& & e_j \cdot \varphi= \mathrm{i}  q_j \varphi,\nonumber\\
& & e_{n+j} \cdot \varphi=\partial_{q_j} \varphi
\end{eqnarray}
for $j=1\ldots,n$.
The three differential operators valued in $\mathrm{End}({\mathcal S}({\mathbb R}^n))$,
\begin{eqnarray}\label{operDsxsE}
& & X_s={} \sum_{j=1}^n (y_{j} \partial_{q_j} + \mathrm{i}  x_j q_j),
\nonumber \\
& & D_s={} \sum_{j=1}^n (\mathrm{i}  q_j \partial_{y_{j}} - \partial_{x_j}\partial_{q_j}),
\nonumber \\
& & E={} \sum_{j=1}^{n} (x_{j}\partial_{x_j}+y_{j}\partial_{y_j}),
\end{eqnarray}
are $\mathfrak{sp}(2n,\mC)$-equivariant and  
generate the representation of the Lie algebra 
$\sl(2)$ on the space $\mathrm{Pol}(\mR^{2n},\mC)\otimes {\mathcal S}(\mR^n)$. 
Their commutation relations are
\begin{eqnarray}\label{comrel}
& & [E +n,D_s]=-D_s,
\nonumber \\
\label{slRels}
& & [E +n,X_s]=X_s,\nonumber \\
& & [X_s,D_s]=\mathrm{i} (E +n).
\end{eqnarray}
The metaplectic analogue of the classical theorem on the separation of variables allows to decompose 
polynomial symplectic spinors $\mathrm{Pol}(\mR^{2n},\mC)\otimes {\mathcal S}(\mR^n)$ under the action 
of $\mp(2n,\mR)$ into a direct sum of simple (irreducible) $\mp(2n,\mR)$-modules, cf. 
\cite{bss}:
\begin{eqnarray}
\mathrm{Pol}(\mR^{2n},\mC)\otimes {\mathcal S}(\mR^n)\simeq\bigoplus_{l=0}^\infty\bigoplus_{j=0}^\infty X_s^j{M}_l^s
\end{eqnarray} 
with
\[
 {M}_l^s:=\big(\mathrm{Pol}_l(\mR^{2n},\mC)\otimes {\mathcal S}(\mR^n)\big)
\cap\mathrm{Ker}(D_s),
\]
which can be graphically represented as
\begin{eqnarray}\label{obrdekonposition}\nonumber
\xymatrix@=11pt{P_0 \otimes {\mathcal S} \ar@{=}[d] &  P_1 \otimes {\mathcal S} \ar@{=}[d]& 
P_2 \otimes {\mathcal S} \ar@{=}[d] & P_3 \otimes {\mathcal S} \ar@{=}[d] & 
P_4 \otimes {\mathcal S} \ar@{=}[d]& P_5 \otimes {\mathcal S} \ar@{=}[d] \\
M_0^s \ar[r] & X_s M_0^s \ar @{} [d] |{\oplus} \ar[r] & X_s^2 M_0^s \ar @{} [d] |{\oplus} \ar[r] & X_s^3 M_0^s \ar @{} [d] |{\oplus}
 \ar[r] & X_s^4 M_0^s \ar @{} [d] |{\oplus}\ar[r] & X_s^5 M_0^s \ar @{} [d] |{\oplus} \\
& M_1^s \ar[r] & X_s M_1^s \ar @{} [d] |{\oplus}\ar[r] & X_s^2 M_1^s \ar @{} [d] |{\oplus}
 \ar[r] & X_s^3 M_1^s \ar @{} [d] |{\oplus}\ar[r] & X_s^4 M_1^s \ar @{} [d] |{\oplus} \\
&& M_2^s \ar[r] & X_s M_2^s \ar @{} [d] |{\oplus}
 \ar[r] & X_s^2 M_2^s \ar @{} [d] |{\oplus}\ar[r] & X_s^3 M_2^s \ar @{} [d] |{\oplus} \\
&&& M_3^s \ar[r] & X_s M_3^s \ar @{} [d] |{\oplus}\ar[r] & X_s^2 M_3^s  \ar @{} [d] |{\oplus} \\
&&&& M_4^s \ar[r] & X_s M_4^s \ar @{} [d] |{\oplus}  \\
&&&&& M_5^s 
}
\end{eqnarray}
To simplify the scheme, we used the notation $P_j$ instead 
of $\mathrm{Pol}_j(\mR^{2n},\mC)$ and $\mathcal{S}$ 
instead of $\mathcal{S}(\mR^n)$ in the last picture.
The symplectic Dirac operator $D_s$ as well as $X_s$ act horizontally in 
the previous picture, but in opposite 
directions; $E$ preserves each simple metaplectic 
module in the decomposition.


\section{Symmetries of the symplectic Dirac operator}

We shall start the present section with a short reminder of the 
notion of symmetry operators for the classical Dirac operator 
associated to a quadratic form, see \cite{Eastwood} and \cite{ESS},
and then pass to the case of our interest: the symplectic Dirac operator.

The Clifford algebra associated to a vector space equipped with 
a quadratic form $B$ is determined by the relations 
$e_j\cdot e_k+e_k\cdot e_j=-2 B(e_j,e_k)$, while the relations for 
the symplectic Clifford algebra 
on $(\mR^{2n},\omega)$ are introduced in Definition \ref{cliffalgdef}. 
In the orthogonal case, the Dirac operator on $\mR^m$  
is $D=\sum_{j=1}^m e_j \partial_{x_j}$ and its polynomial solutions are coined spherical 
monogenics. The module of polynomial spherical 
monogenics of homogeneity $h$ is denoted by 
$M_h=\big(\mathrm{Pol}_h(\mR^m,\mC) \otimes \mathbb{S}\big)\cap\mathrm{Ker}(D)$, where 
$\mS$ is the spinor space.
In particular, each of the modules $M_h$, $h\in \mN_0$, is an irreducible representation of the 
Lie algebra $\mathfrak{so}(m)$ acting by the differential
operators $$K_{jk}=x_j\partial_{x_k}-x_k \partial_{x_j}-\frac{1}{2}e_j e_k, \quad j\neq k,\quad j,k=1,\ldots,m .$$
Moreover, the space $M=\bigoplus_{h} M_h$ is an irreducible 
representation of the conformal Lie algebra $\mathfrak{so}(m+1,1,\mR)$, which is
the linear span of 
$K_{jk}$, $2E+m-1$, $\partial_{x_j}$ and $\tilde{T}_j$, $j,k=1,\ldots,m$;
here the operators $\tilde{T}_j: M_h\rightarrow M_{h+1}$ act by 
\begin{equation}
\tilde{T}_j=X e_j+x_j(m+2E)-|X|^2 \partial_{x_j},
\end{equation}
where
\[
X=\sum_{j=1}^m e_j x_j.
\]

Let us now turn our attention to the symplectic space $(\mR^{2n},\omega)$. 
First of all, we find differential operators increasing  
the homogeneity 
of polynomial solutions of the symplectic Dirac operator by one. We construct them
as a composition of the multiplication by $x_l, y_l$, $l=1,\ldots,n$, and projection on 
the kernel of the symplectic Dirac operator $D_s$.

In the first step we check that $D_s^3$ acts trivially on $x_l m$, $y_l m$
for $m\in M_{h}^s$ and coordinate functions $x_l$, $y_l$ on $\mR^{2n}$. 
For $j=1,\ldots,n$, we have  
\begin{eqnarray}\label{Ds2xf}
&&D_s^2(x_j m)=D_s(-\partial_{q_j}m)=-\sum_{k=1}^{2n}(\mathrm{i}  q_k \partial_{y_{k}}
-\partial_{q_k}\partial_{x_k})\partial_{q_j}m=\mathrm{i}  \partial_{y_{j}}m, \nonumber\\
&&D_s^2(y_{j} m)=D_s(\mathrm{i}q_jm)=
\sum_{k=1}^{2n}(\mathrm{i} q_k \partial_{y_{k}}-\partial_{q_k}\partial_{x_k})\mathrm{i}  q_jm
=-\mathrm{i}  \partial_{x_{j}}m,
\end{eqnarray}
and so $x_l m$, $y_l m$ are in the kernel of $D_s^3$ for all $l=1,\ldots ,n$. 
Denoting the identity endomorphism $\Id$, the 
corresponding projector of $x_l m$, $y_l m$ on the 
homogeneity $h+1$ subspace of $\mathrm{Ker}(D_s)$ is 
$$P^s_{h+1}=\Id+cX_sD_s+dX_s^2 D_s^2$$ 
for some constants $c, d$ depending on $h$ and $n$. 
The relations \eqref{slRels} imply that
on the spaces of homogeneous symplectic monogenics 
holds
\begin{eqnarray}\label{symplproj}
 P_{h+1}^s m_{h+1} &=& m_{h+1},\nonumber \\
 P_{h+1}^s X_s m_{h} &=& X_s m_h +cX_sD_s X_sm_h=\big(1-\mathrm{i}  c(h+n)\big) X_s m_h,\nonumber\\
 P_{h+1}^s X_s^2 m_{h-1} &=& X_s^2 m_{h-1}
 +cX_s D_s(X_s^2 m_{h-1})+dX_s^2 D_s^2 (X_s^2m_{h-1})\nonumber\\
 &=& X_s^2m_{h-1}-\mathrm{i}  cX_s^2(h-1+n)m_{h-1}-\mathrm{i} cX_s(h+n)X_sm_{h-1}\nonumber\\
 & &  - d (2h+2n-1)(h+n-1)X_s^2m_{h-1}. 
\end{eqnarray}
Then the second and the third expressions in \eqref{symplproj} 
are zero provided
$$c=\frac{1}{\mathrm{i} (h+n)},\quad d=\frac{-1}{(h+n)(2h+2n-1)},$$ 
hence the projector is
\begin{equation}
P_{h+1}^s=\Id+\frac{1}{\mathrm{i} (h+n)}X_sD_s-\frac{1}{(h+n)(2h+2n-1)}X_s^2 D_s^2.
\end{equation}
The action of the operators $S_l=P_{h+1}^sx_l$, $l=1,\dots ,n$ and $S_{n+l}=P_{h+1}^sy_l$, $l=1,\dots ,n$, on $m\in M_h^s$ is then
\begin{eqnarray*}
S_j m &=& x_j m-cX_s \partial_{q_j}m+ \mathrm{i}\,  d X_s^2 \partial_{y_{j}}m,\\
S_{n+j} m &=& y_{j} m+c X_s \mathrm{i}q_jm- \mathrm{i}\,  d X_s^2 \partial_{x_{j}}m,
\end{eqnarray*}
so we can define for $j=1,\dots ,n$ the collection of differential operators
\begin{eqnarray}
 Z_j &:=& -\mathrm{i} (h+n)(2h+2n-1)S_{n+j},\nonumber\\
 Z_{n+j} &:=& \mathrm{i} (h+n)(2h+2n-1)S_{j}.
\end{eqnarray}

\begin{proposition}\label{ZinDimN}
Let $n\in\mN$. The differential operators
\begin{eqnarray}
 Z_j &=& X_s^2\partial_{x_j}-\mathrm{i}  y_{j} (E +n)(2E+2n-1)-\mathrm{i} X_s q_j (2E+2n-1),
\nonumber \\
 Z_{n+j} &=& X_s^2\partial_{y_{j}}+\mathrm{i}  x_j (E +n)(2E+2n-1)-X_s \partial_{q_j} (2E+2n-1)
\nonumber \\
\end{eqnarray}
for $j=1,\ldots ,n$ are $\mp(2n,\mR)$-equivariant and
preserve the solution space of the symplectic Dirac operator on $(\mR^{2n},\omega)$. 
The operators $Z_l$, $l=1,\ldots, 2n$, increase the homogeneity in the basis variables 
$x_1,\dots ,x_{n}$, $y_1,\dots ,y_{n}$ by one: 
\begin{eqnarray}
& & Z_l: \mathrm{Ker}{(D_s)} \to \mathrm{Ker}(D_s),
\nonumber \\
& & Z_l: M_h\mapsto M_{h+1}, \quad l=1,\dots ,2n,
\end{eqnarray}
where $M_h$ is the irreducible $\mp(2n,\mR)$-module of homogeneity $h$ symplectic polynomial 
spinors in $\mathrm{Ker}(D_s)$.
\end{proposition}
{\bf Proof:}
The property of $\mp(2n,\mR)$-equivariance means that the vector space of dimension 
$2n$ generated by $\{Z_1,\ldots , Z_n, Z_{n+1},\ldots , Z_{2n}\}$ transforms in the
fundamental vector representation of $\mp(2n,\mR)$ with respect to the canonical Lie
algebra structure on the associative algebra $W_{4n}\otimes Cl_s(\mR^{2n},\omega)$.
Recall that $W_{4n}$ is generated by $x_j, y_j, \partial_{x_{j}}, \partial_{y_{j}}$ 
for $j=1, \ldots ,n$. The verification of all commutation relations of the Lie algebra
$\mp(2n,\mR)$ (cf., \eqref{ssrepr}) with $\{Z_1,\ldots , Z_n, Z_{n+1},\ldots , Z_{2n}\}$
is a straightforward but tedious computation.

The second part of the claim is a consequence of   
\begin{eqnarray}
\, [D_s, X^2_s\partial_{x_j}] &=& -\mathrm{i}  X_s \partial_{x_j}(2E+2n-1), \nonumber \\
\, [D_s,\omega_{jk}\delta^{k,l}x_l(E+n)(2E+2n-1)] &=& -e_j(E+n)(2E+2n-1)
\nonumber \\
& & +\omega_{jk}\delta^{k,l}x_l(4E+4n+1)D_s,
\nonumber \\
\, [D_s,X_s e_j(2E+2n-1)] &=& -\mathrm{i}  e_j(E+n)(2E+2n-1)\nonumber \\
& & -\mathrm{i}  X_s \partial_{x_j}(2E+2n-1)+2X_s e_j D_s, 
\end{eqnarray}
because the linear combination 
$$
AX^2_s\partial_{x_j} +B\omega_{jk}\delta^{k,l}x_l(E+n)(2E+2n-1) +CX_se_j(2E+2n-1) 
$$
for $A,B,C\in\mC$ and all $j=1,\dots ,2n$ commutes with $D_s$ provided 
$A=1, B=\mathrm{i} $ and $C=-1$. 
To shorten our notation we used $\omega_{jk}=\omega(e_j,e_k)$, see \eqref{omega}.

\hfill
$\square$

The differential operators $Z_j$, $Z_{n+j}$, $j=1,\ldots,n$ are of third order, 
and are of second order in the base variables $x_j$, $y_j$ (due to their quadratic 
dependence on the homogeneity operator $E$.)

\begin{proposition}\label{DxjinDimN}
The $\mp(2n,\mR)$-equivariant first order differential operators 
\begin{eqnarray}
\partial_{x_j}, \partial_{y_j},\quad j=1,\dots ,n 
\end{eqnarray}
preserve the solution space of the symplectic Dirac operator on $(\mR^{2n},\omega)$: 
\begin{eqnarray}
& & \partial_{x_j},\partial_{y_j}: \mathrm{Ker}(D_s) \to \mathrm{Ker}(D_s),\nonumber \\
& & \partial_{x_j},\partial_{y_j}: M_h\mapsto M_{h-1}, \quad j=1,\dots ,n,
\end{eqnarray}
where $M_h$ is the irreducible $\mp(2n,\mR)$-module of homogeneity $h$ symplectic polynomial 
spinors in $\mathrm{Ker}(D_s)$.
\end{proposition}
{\bf Proof:}
The property of $\mp(2n,\mR)$-equivariance means that the vector space of dimension 
$2n$ generated by $\{\partial_{x_1},\ldots ,\partial_{x_n}, \partial_{y_1}, \ldots,
\partial_{y_n}\}$ transforms in the
fundamental vector representation of $\mp(2n,\mR)$ with respect to the canonical Lie
algebra structure on the associative algebra $W_{4n}\otimes Cl_s(\mR^{2n},\omega)$.
The verification of all commutation relations of the Lie algebra $\mp(2n,\mR)$ 
(cf., \eqref{ssrepr}) with $\{\partial_{x_1},\ldots ,\partial_{x_n}, \partial_{y_1}, \ldots,
\partial_{y_n}\}$ is a straightforward computation.

The rest of the claim follows from $[\partial_{x_j},D_s]=0$ and $[\partial_{y_j},D_s]=0$ 
for $j=1,\dots ,n$.

\hfill
$\square$

\subsection{First order symmetries of the symplectic Dirac operator on $(\mR^{2},\omega)$}

The aim of the present section is to compute all first order differential
operators (in both the horizontal variables $x,y$ and the vertical variable $q$) 
which are symmetries of the symplectic Dirac operator. Here we restrict to $n=2$, 
the case of general even dimension being notationally tedious.  

We start with $(\mR^{2},\omega)$ and denote the coordinates 
by $x=x_1, y=y_1$, the coordinate vector fields by 
$\partial_{x},\partial_{y}$ and a symplectic frame is $e_{1},e_{2}$ 
with the action on a symplectic 
spinor $\varphi\in \mathrm{Pol}(\mR^2,\mC)\otimes \mathcal{S}(\mR)$ 
\begin{eqnarray*}
e_1 \cdot\varphi= \mathrm{i}  q \varphi,\quad e_{2} \cdot\varphi=\partial_q \varphi.
\end{eqnarray*}
Following \eqref{ssrepr}, the basis elements of 
$\mp(2,\mR)(\simeq \sp(2,\mR)\simeq {\mathfrak sl}(2)$) act as
\begin{eqnarray}\label{HXYxy}
& & \tilde{X}=- y \partial_x-\frac{\mathrm{i} }{2} q^2,\,\nonumber \\
& &\tilde{Y}=- x \partial_y-\frac{\mathrm{i} }{2} \partial_q^2,\,
\nonumber \\
& & \tilde{H}=- x \partial_x+y \partial_y+q \partial_q+\frac{1}{2},
\end{eqnarray}
and satisfy the commutation relations of the Lie algebra 
$\mp(2,\mR)$:
\begin{align*}
[\tilde{X},\tilde{Y}]&=\tilde{H},\\
[\tilde{H},\tilde{X}]&=2\tilde{X},\\
[\tilde{H},\tilde{Y}]&=-2\tilde{Y}.
\end{align*}
Notice that these operators preserve homogeneity in the variables $x,y$.
The three $\mp(2,\mR)$-invariant operators 
\begin{eqnarray}
&& X_s = y \partial_q + \mathrm{i}  x q,
\nonumber \\
&& D_s = \mathrm{i}  q \partial_y - \partial_x\partial_q,
\nonumber \\
&& E  = x\partial_x+y\partial_y
\end{eqnarray}
form the Lie algebra isomorphic to $\sl(2)$.
The operators $X_s, D_s$ and $E$ commute with  
$\tilde{X},\tilde{Y}$ and $\tilde{H}$, i.e.
they are $\mp(2,\mR)$ intertwining differential 
operators on complex polynomials valued in the 
Segal-Shale-Weil representation. 
A consequence of Proposition \ref{ZinDimN} 
and Proposition \ref{DxjinDimN} is
\begin{corollary}\label{Z1Z2xy}
The commuting operators 
\begin{eqnarray}
& & Z_1=-X_s^2\partial_x+ \mathrm{i} y (E +1)(2E+1)+X_s \mathrm{i}q (2E+1), \nonumber \\
& & Z_2=-X_s^2\partial_y- \mathrm{i} x (E +1)(2E+1)+X_s \partial_q (2E+1)
\end{eqnarray}
preserve the solution space of the symplectic Dirac operator $D_s$ and increase 
the homogeneity in the variables $x,y$ by one, $Z_j:M_h\mapsto M_{h+1}$, 
$j=1,2$, for $M_h$ being the irreducible $\mp(2,\mR)$-module of homogeneity $h$ 
polynomial symplectic spinors in $\mathrm{Ker}(D_s)$.

The commuting operators 
\begin{equation}
\partial_x, \partial_y
\end{equation}
preserve the solution space  of the symplectic Dirac operator $D_s$ and decrease the homogeneity 
in the variables $x,y$ by one.

\end{corollary}

The commutator $[Z_1,Z_2]$ is zero, and
\begin{eqnarray}
& &  [\partial_x,Z_1] 
=-2\mathrm{i} \tilde{X}(2E+1), \nonumber \\
& & [\partial_y,Z_1] = 2X_sD_s+\mathrm{i} \tilde{H}(2E+1)+
\mathrm{i}  (2E+1)(2E+1)+\frac{\mathrm{i}}{2},\nonumber \\
& &  [\partial_x,Z_2] = -2X_sD_s+\mathrm{i} \tilde{H}(2E+1)
-\mathrm{i}  (2E+1)(2E+1)-\frac{\mathrm{i} }{2},
\nonumber\\
& & [\partial_y,Z_2]  =2\mathrm{i} \tilde{Y}(2E+1).
\end{eqnarray}
Moreover, we have
\begin{equation}
\begin{array}{ll}
\, [Z_1,\tilde{H}]=-Z_1, & [Z_2,\tilde{H}]=Z_2, \\
\, [Z_1,\tilde{X}]=0,    & [Z_2,\tilde{X}]=-Z_1, \\
\, [Z_1,\tilde{Y}]=Z_2,  & [Z_2,\tilde{Y}]=0,  \\
\, [Z_1,E]=-Z_1,         & [Z_2,E]=-Z_2,
\end{array}
\end{equation}
as well as
\begin{equation}
\begin{array}{ll}
\, [\partial_x,\tilde{H}]=-\partial_x  & [\partial_y,\tilde{H}]=\partial_y,  \\
\, [\partial_x,\tilde{X}]= 0   & [\partial_y,\tilde{X}]=-\partial_x, \\
\, [\partial_x,\tilde{Y}]=-\partial_y  & [\partial_y,\tilde{Y}]=0,  \\
\, [\partial_x,E]=\partial_x,  & [\partial_y,E]=\partial_y.  
\end{array}
\end{equation}

\begin{remark}
The commutator of commutators 
$[\partial_x,Z_1]=-2\mathrm{i} \tilde{X}(2E+1)$ and $[\partial_y,Z_2]=2\mathrm{i} \tilde{Y}(2E+1)$ 
gives
$$[-2\mathrm{i} \tilde{X}(2E+1),2\mathrm{i} \tilde{Y}(2E+1)]=4\tilde{H}(2E+1)(2E+1).$$
Then we can compute the commutator of this commutator with, for example,
$[\partial_x,Z_1]=-2\mathrm{i} \tilde{X}(2E+1)$, resulting in the third power of $(2E+1)$. 
In general, we can produce an arbitrarily high power of $(2E+1)$ in iterated commutators, 
hence the linear span of the operators 
$\tilde{H},\tilde{X}, \tilde{Y}, \partial_x, \partial_y, Z_1, Z_2$ and $E$ 
is not closed under the commutator bracket. 
\end{remark}

Let us briefly mention the key concept of (generalized) differential symmetries for 
the symplectic Dirac operator, see \cite{Eastwood} and references therein for 
an introduction.  
A differential operator 
$A$ is a symmetry of $D_s$ if there exists another differential operator $B$ 
such that 
\begin{eqnarray}\label{gensymdef}
D_s A=B D_s.
\end{eqnarray}
Consequently, symmetry operators preserve the solution space of the symplectic Dirac operator.
\begin{theorem}
The first order symmetries of the symplectic Dirac operator $D_s$ 
on $\mR^{2}$ are given by the linear span of differential operators 
$\partial_x$, $\partial_y$, 
$\tilde{H}$, $\tilde{X}$, $E$ and 
$y\tilde{H}-2x \tilde{X}+yE+\frac{3}{2}y$.
\label{Symm_sympl_thm}
\end{theorem}
{\bf Proof:}
Let us consider a general first order differential operator in the variables $x,y,q$: 
$$
A=F_0(x,y,q)\partial_x+F_1(x,y,q)\partial_y+F_2(x,y,q)\partial_q+F_3(x,y,q),
$$
where $F_j, j=0,1,2,3$, are convenient functions of $x,y$ and $q$. 
Then $D_sA=A D_s+[D_s,A]$, so that \eqref{gensymdef} implies 
$[D_s,A]=B' D_s$ for a differential operator $B'$. 
The computation of commutators gives
\begin{eqnarray*}
&& \big(\mathrm{i}q[\partial_y,F_0(x,y,q)]-\partial_q[\partial_x,F_0(x,y,q)]
-[\partial_q,F_2(x,y,q)]\partial_q-[\partial_q,F_3(x,y,q) ] \big)\partial_x\\
&& +\big(\mathrm{i}q[\partial_y,F_1(x,y,q)]-\partial_q [\partial_x,F_1(x,y,q)]+F_2(x,y,q)
[\mathrm{i}q,\partial_q] \big)\partial_y\\
&& - [\partial_q,F_0(x,y,q)]\partial_x^2 - [\partial_q,F_1(x,y,q)]\partial_x\partial_y\\
&& +\mathrm{i}q[\partial_y,F_2(x,y,q)]\partial_q - \partial_q [\partial_x,F_2(x,y,q)]\partial_q + 
\mathrm{i}q[\partial_y,F_3(x,y,q)]-\partial_q[\partial_x,F_3(x,y,q)]\\
&&=B'(\mathrm{i}q\partial_y-\partial_x\partial_q).
\end{eqnarray*}
The commutator $[\partial_q,F_0(x,y,q)]$ by $\partial_x^2$ does not
depend on $\partial_q$ and so equals to zero. Hence $F_0(x,y,q)$ is independent of the variable $q$, $F_0\equiv F_0(x,y)$. Then the 
commutator $[\partial_q,F_1(x,y,q)]$ by $\partial_x\partial_y$ has 
to be zero as well, i.e., $F_1(x,y)$ is independent of $q$. Moreover, 
the commutator $[\partial_x,F_1(x,y)]$ in 
$\partial_q [\partial_x,F_1(x,y)]\partial_y$ has to be zero, i.e., 
$F_1\equiv F_1(y)$.

We can separate the last equation into three equalities:
\begin{align}
& \big(\mathrm{i}q[\partial_y,F_0(x,y)]-[\partial_q,F_3(x,y,q) ]-\big([\partial_x,F_0(x,y)]+[\partial_q,F_2(x,y,q)]\big)\partial_q \big)\partial_x
\nonumber \\
& =-B'\partial_q\partial_x \label{rceII},\\
& \big(\mathrm{i}q[\partial_y,F_1(y)]-\mathrm{i} F_2(x,y,q) \big) \partial_y=B' \mathrm{i}q\partial_y 
\label{rceI},\\
& \mathrm{i}q[\partial_y,F_2(x,y,q)]\partial_q - \partial_q [\partial_x,F_2(x,y,q)]\partial_q + 
\mathrm{i}q[\partial_y,F_3(x,y,q)]
\nonumber \\
& -\partial_q[\partial_x,F_3(x,y,q)]=0. \label{rceIII}
\end{align}
The equation \eqref{rceII} yields $\mathrm{i}q[\partial_y,F_0(x,y)]-[\partial_q,F_3(x,y,q) ]=0$.
We set 
\begin{equation}
F_3(x,y,q)=F_3'(x,y)\frac{\mathrm{i}}{2}q^2+F_3''(x,y), \label{vyjF3}
\end{equation}
and therefore
\begin{eqnarray}
[\partial_y,F_0(x,y)]=F_3'(x,y)\label{rceVI}.
\end{eqnarray}
The second equality \eqref{rceI} implies 
\begin{equation}
F_2(x,y,q)=F_2'(x,y)q. \label{vyjF2}
\end{equation}
Then $[\partial_q,F_2(x,y,q)]=F_2'(x,y)$, and equations \eqref{rceI} and \eqref{rceII} give
\begin{eqnarray}
&& [\partial_y,F_1(y)]-F_2'(x,y)= [\partial_x,F_0(x,y)]+F_2'(x,y),\nonumber\\
&& [\partial_y,F_1(y)]= [\partial_x,F_0(x,y)]+2F_2'(x,y).\label{rceV}
\end{eqnarray}
The equation \eqref{rceIII} can be rewritten with the use of \eqref{vyjF3} and \eqref{vyjF2} as
\begin{eqnarray*}
&& [\partial_y,F_2'(x,y)]\mathrm{i}q^2\partial_q - [\partial_x,F_2'(x,y)](\partial_q+q\partial_q^2)-[\partial_y,F_3'(x,y)]\frac{1}{2}q^3\\
&& +[\partial_y,F_3''(x,y)]\mathrm{i}q - [\partial_x,F_3'(x,y)](\mathrm{i}q +\frac{1}{2}\mathrm{i}q^2\partial_q)-[\partial_x,F_3''(x,y)]\partial_q=0.
\end{eqnarray*}
Because there is only one commutator by $q\partial_q^2$ and $q^3$,
we have $F_2'\equiv F_2'(y)$, $F_3'\equiv F_3'(x)$. Then the commutators  
by $\partial_q$ have to be zero and $F_3''$ is independent of $x$,
$F_3''\equiv F_3''(y)$.
 The commutators by $\mathrm{i}q^2\partial_q$ and $\mathrm{i}q$ give the relations
\begin{equation}\label{rceVV}
[\partial_y,F_2'(y)]-\frac{1}{2}[\partial_x,F_3'(x)]=0, 
\quad \quad [\partial_y,F_3''(y)]-[\partial_x,F_3'(x)]=0.
\end{equation}
The solution of \eqref{rceVV} is $F_2'(y)=\frac{1}{2}\alpha y+\gamma$, 
$F_3'(x)=\alpha x+\beta$ and $F_3''=\alpha y+\gamma$. The substitution of 
this solution into \eqref{rceVI} yields $F_0(x,y)=\alpha xy+ \beta y+F_0'(x)$. 
Substituting into \eqref{rceV}, we get $F_0'(x)=\eta x+ \zeta$ and 
$F_1(y)=\alpha y^2 +(2\gamma+\eta)y+\kappa$.
Taken altogether, the functions $F_j, j=0,1,2,3$ are
\begin{eqnarray*}
F_0=\alpha xy +\eta x+\beta y+\zeta, \quad && F_2=\frac{1}{2}\alpha yq+\gamma q, \\
F_1=\alpha y^2 +(2\gamma+\eta)y+\kappa, \quad && 
F_3=\left(\alpha x+\beta \right)\frac{\mathrm{i}}{2}q^2 +\alpha y+\delta,
\end{eqnarray*}
where $\alpha,\beta,\gamma,\delta,\eta,\zeta,\kappa\in \mathbb{C}$ are arbitrary constants. 
The constant $\beta$ corresponds to the operator $\tilde{X}$, $\zeta$ and 
$\kappa$ correspond to $\partial_x$ and $\partial_y$. A combination of $\eta, \gamma$ and 
$\delta$ corresponds to a combination of $E$, $\tilde{H}$ and the identity operator. Finally, $\alpha$ corresponds to the operator $y\tilde{H}-2x \tilde{X}+yE+\frac{3}{2}y$.

\hfill
$\square$

We notice that $\tilde{Y}$ is a second order differential operator, but it is
first order in the base variables $x,y$. 
The operators 
$Z_1,Z_2$ are symmetries of $D_s$ but they are third order differential operators, 
second order in 
the base variables $x,y$. 


\subsection{First order symmetries in 
the holomorphic variable}

\label{opervzbarz}
We use the complex coordinates $z=x+\mathrm{i} y,\overline{z}=x-\mathrm{i} y$,
for the standard complex structure on $\mR^2$, where 
$\partial_x=\partial_z+\partial_{\bar{z}}$ and 
$\partial_y=\mathrm{i} (\partial_z-\partial_{\bar{z}})$. 
In the complex coordinates $z, \bar{z}$ we have 
\begin{eqnarray}
& & X_s=\frac{\mathrm{i}}{2}\big( (q-\partial_q)z+(q+\partial_q)\bar{z}\big),\label{XsHolom}
\nonumber \\
& & D_s=-(q+\partial_q)\partial_z+(q -\partial_q)\partial_{\bar{z}} ,
\label{DsHolom} 
\nonumber \\
& & E= z\partial_z+ \bar{z}\partial_{\bar{z}} 
\end{eqnarray}
and 
\begin{eqnarray}\label{Z1Z2zbarz}
&&Z_1=2X_s^2\partial_z+\bar{z}(E+1)(2E+1)+\mathrm{i} X_s(\partial_q-q)(2E+1),\nonumber\\
&&Z_2=2X_s^2\partial_{\bar{z}}-z(E+1)(2E+1)-\mathrm{i} X_s(\partial_q+q)(2E+1),\quad
\end{eqnarray}
 where $Z_1=\bar{Z_1 }+\mathrm{i} \bar{Z_2}$ and $Z_2=\bar{Z_1}-\mathrm{i} \bar{Z_2}$, cf. 
$\bar{Z_1}$ and $\bar{Z_2}$ in Corollary \ref{Z1Z2xy}.

The commutator of $[Z_1,Z_2]$ is trivial and the commutators with (anti-)holomorphic
coordinate vector fields are
\begin{eqnarray}
&& [Z_1,\partial_{z}]=2\mathrm{i} X_t (2E+1),\nonumber \\
&& [Z_1,\partial_{\bar{z}}]=2 \mathrm{i} X_s D_s -H_t (2E+1)-(2E+1)(2E+1)-\frac{1}{2},\nonumber \\
&& [Z_2,\partial_{z}]=-2\mathrm{i} X_s D_s -H_t (2E+1)+(2E+1)(2E+1)+\frac{1}{2},\nonumber \\
&& [Z_2,\partial_{\bar{z}}]=2\mathrm{i} Y_t (2E+1), 
\end{eqnarray}
where we introduced 
\begin{eqnarray*}
&&H_t=\mathrm{i}\bar{\tilde{X}}-\mathrm{i}\bar{\tilde{Y}},\\
&&X_t=-\frac{1}{2}\big(\bar{\tilde{X}}+\bar{\tilde{Y}} +\mathrm{i} \bar{\tilde{H}}\big),\\
&&Y_t=-\frac{1}{2}\big(\bar{\tilde{X}}+\bar{\tilde{Y}} -\mathrm{i} \bar{\tilde{H}}\big),
\end{eqnarray*}
with $\bar{\tilde{H}},\bar{\tilde{X}}$ and $ \bar{\tilde{Y}}$ the 
operators \eqref{HXYxy} in the variables $z,\bar{z}$:
\begin{eqnarray}\label{HXYholom}
&& H_t=\bar{z}\partial_{\bar{z}}-z\partial_{z}+\frac{1}{2}(q^2-\partial_q^2),\nonumber\\
&& X_t=\mathrm{i}\bar{z}\partial_{z}+\frac{\mathrm{i}}{4}(q-\partial_q)^2,\nonumber\\
&& Y_t=-\mathrm{i} z\partial_{\bar{z}}+\frac{\mathrm{i}}{4}(q+\partial_q)^2.
\end{eqnarray}
The operators $H_t,X_t$ and $Y_t$ commute with $D_s, X_s, E$, and satisfy 
the commutation relations of the Lie algebra $\mp(2,\mR)$:
\begin{eqnarray}
&&[X_t,Y_t]=H_t,\nonumber\\
&&[H_t,X_t]=2X_t,\nonumber\\
&&[H_t,Y_t]=-2Y_t.
\end{eqnarray}
A straightforward computation reveals 
\begin{equation}
\begin{array}{ll}
\, [Z_1,H_t]=-Z_1, &[Z_2,H_t]=Z_2, \\
\, [Z_1,X_t]=0,    &[Z_2,X_t]=\mathrm{i} Z_1, \\ 
\, [Z_1,Y_t]=-\mathrm{i} Z_2, &[Z_2,Y_t]=0, \\
\, [Z_1,E]=-Z_1,   &[Z_2,E]=-Z_2, \\
\end{array}
\end{equation}
\begin{equation}
\begin{array}{ll}
\, [\partial_z,H_t]=-\partial_z,			& [\partial_{\bar{z}},H_t]=\partial_{\bar{z}}, \\
\, [\partial_z,X_t]=0,   					&[\partial_{\bar{z}},X_t]=\mathrm{i} \partial_z, \\ 
\, [\partial_z,Y_t]=-\mathrm{i} \partial_{\bar{z}},	&[\partial_{\bar{z}},Y_t]=0, \\
\, [\partial_z,E]=\partial_{z},	            &[\partial_{\bar{z}},E]=\partial_{\bar{z}}. \\
\end{array}
\end{equation}


\section{Towards a symplectic Clifford-Fourier transform }
\label{towarsdSFT}

The central role in harmonic analysis on $\mR^n$ is played by the 
Lie algebra $\mathfrak{sl}(2,\mC)$, generated by the $\mathfrak{so}(n,\mR)$-invariant 
Laplace operator $\triangle$ and the norm squared $|x|^2$ of the 
vector $x\in\mathbb{R}^n$. 
The classical integral Fourier transform, 
\begin{eqnarray}
F(f)(y) = (2\pi)^{-\frac{n}{2}}\int_{\mathbb{R}^n}f(x)
\exp^{-\mathrm{i}\langle x,y\rangle } \,\mathrm{d}x, \quad \langle x,y\rangle =
\sum\limits_{i=1}^n x_iy_i, 
\end{eqnarray}
can be equivalently represented by the operator exponential that contains the generators
of $\mathfrak{sl}(2,\mC)$:
\begin{eqnarray}
\exp^{\frac{\mathrm{i}\pi n}{4}}\exp^{\frac{\mathrm{i}\pi}{4}(\triangle -|x|^2)},
\end{eqnarray}
which means that the two operators have the same spectral properties.
There are analogous results in the harmonic analysis for finite groups based on Dunkl operators,
or Clifford analysis based on the Clifford algebra associated to a quadratic form and the Dirac
operator $D=\sum\limits_{j=1}^n e_j\partial_{x_j}$, written in a basis $e_1,\ldots e_n$ of 
$\mathbb{R}^n$ with coordinates $x_1,\ldots, x_n$,
cf. \cite{DOSSDunkl}, \cite{DOSSdef} and \cite{DBX}.

In the present section we discuss several basic questions in this direction, focusing
on symplectic Clifford analysis and the associated symplectic Dirac operator in real 
dimension $2$. 

\subsection{The eigenfunction decomposition for the operator $D_s - c X_s$}

The symplectic Fourier transform is based on the eigenvalue equation
\begin{equation}\label{eigensf}
(D_s - c X_s)f=\lambda f,\quad c\in \mR, \lambda\in\mC .
\end{equation}
As already indicated, we shall stick to the real dimension $2$ and 
look for the solutions of this equation in terms of a linear 
combination of elements $g(X_s)m_k^s$, where $m_k^s\in M_k^s$ is a symplectic monogenic
and $g$ is a polynomial in the variable $X_s$. 
We shall first focus on the problem whether for a symplectic spinor 
$\varphi$ valued in 
$\mathcal{S}(\mR)$ holds $e^{\alpha X_s}\varphi \in \mathcal{S}(\mR)$ for 
$\alpha\in \mC$.
\begin{lemma}\label{Xse^minusqpullemma}
The following identity holds,
\begin{equation}
e^{\alpha X_s}e^{-\frac{q^2}{2}}=e^{-\frac{q^2}{2}} e^{\frac{1}{2}\alpha(\mathrm{i} x-y)(2q+\alpha y)}.
\end{equation}
\end{lemma}
{\bf Proof:}
Writing the exponential as
$$e^{\alpha X_s}e^{-\frac{q^2}{2}}=\sum_{k=0}^\infty \frac{\alpha^k}{k!} X_s^k e^{-\frac{q^2}{2}},$$
we show by induction on $k\in\mN_0$ that 
\begin{equation}\label{Xs^keminusqpul}
X_s^k e^{-\frac{q^2}{2}} =e^{-\frac{q^2}{2}} \sum_{m=0}^{\lfloor \frac{k}{2}\rfloor} 
\frac{k! (\mathrm{i} x - y)^{k-m}q^{k-2m}y^m}{m! (k-2m)!2^m}.
\end{equation}
Recall the notation $\lfloor \cdot \rfloor$ for the floor function. 
The equation is satisfied for $k=0$ and for $k=1$, $X_s e^{-\frac{q^2}{2}}$ 
is equal to $e^{-\frac{q^2}{2}} q(\mathrm{i}x-y)$. 
Assuming \eqref{Xs^keminusqpul} holds for $k$, we aim to prove the identity for $k+1$. 
Let us start with odd $k$:
\begin{eqnarray*}
&& (\mathrm{i}x q + y\partial_q)e^{-\frac{q^2}{2}} \sum_{m=0}^{\lfloor 
\frac{k}{2}\rfloor} \frac{k! (\mathrm{i} x - y)^{k-m}q^{k-2m}y^m}{m! (k-2m)!2^m}\\
&& = e^{-\frac{q^2}{2}} \sum_{m=0}^{\lfloor \frac{k}{2}\rfloor} 
\frac{k! (\mathrm{i} x - y)^{k+1-m}q^{k+1-2m}y^m}{m! (k-2m)!2^m}\\
&& + e^{-\frac{q^2}{2}} \sum_{m=0}^{\lfloor \frac{k}{2}\rfloor} 
\frac{k! (\mathrm{i} x - y)^{k-m}q^{k-1-2m}y^{m+1}}{m! (k-2m-1)!2^m},
\end{eqnarray*}
and the shift $m\mapsto m-1$ in the second sum results into
\begin{eqnarray*}
&& e^{-\frac{q^2}{2}} \sum_{m=0}^{\lfloor \frac{k}{2}\rfloor} 
\frac{(k+1)! (\mathrm{i} x - y)^{k+1-m}q^{k+1-2m}y^m}{m! (k+1-2m)!2^m}
\left(\frac{k+1-2m}{k+1}+\frac{2m}{k+1}\right)\\
&& + e^{-\frac{q^2}{2}} \frac{(k+1)! 
(\mathrm{i} x - y)^{k+1-\frac{k+1}{2}}y^{\frac{k+1}{2}}}
{\left(\frac{k+1}{2}\right)! 2^{\frac{k+1}{2}}} \\
&& =e^{-\frac{q^2}{2}} \sum_{m=0}^{\lfloor \frac{k+1}{2}\rfloor} 
\frac{(k+1)! (\mathrm{i} x - y)^{k+1-m}q^{k+1-2m}y^m}{m! (k+1-2m)!2^m}
\end{eqnarray*}
which proves the induction step. For $k$ even, 
$\lfloor \frac{k}{2}\rfloor=\lfloor \frac{k+1}{2}\rfloor$ 
and the second expression on the last display is zero, so that
\begin{equation*}
e^{\alpha X_s}e^{-\frac{q^2}{2}}=e^{-\frac{q^2}{2}} \sum_{k=0}^\infty \sum_{m=0}^{\lfloor \frac{k}{2}\rfloor} \frac{\alpha^k (\mathrm{i} x - y)^{k-m}q^{k-2m}y^m}{m! (k-2m)!2^m}.
\end{equation*}
The change of the order in the last summation while keeping $m$ fixed gives
\begin{equation}
\sum_{k=2m}^\infty  \frac{\alpha^k (\mathrm{i} x - y)^{k-m}q^{k-2m}y^m}{m! (k-2m)!2^m}= 
\frac{\alpha^{2m}(\mathrm{i}x-y)^my^m}{m! 2^m} e^{\alpha q (\mathrm{i}x-y)},
\end{equation}
and so
\begin{equation*}
e^{-\frac{q^2}{2}} \sum_{m=0}^\infty \frac{\alpha^{2m}(\mathrm{i}x-y)^my^m}{m! 2^m} 
e^{\alpha q (\mathrm{i}x-y)}= e^{-\frac{q^2}{2}} e^{\frac{1}{2}\alpha(\mathrm{i} x-y)(2q+\alpha y)}.
\end{equation*}

\hfill
$\square$

By Lemma \ref{Xse^minusqpullemma} we see that $e^{\alpha X_s}\varphi$, $\alpha\in \mC$,
is for $\varphi=e^{-\frac{q^2}{2}}$ a Schwartz function in the variable 
$q$ and a non-polynomial function in the variables $x,y$. 
This property remains true for any $\varphi=p(x,y) e^{-\frac{q^2}{2}}$, where 
$p(x,y)\in \mathrm{Pol}(\mR^2,\mC)$:
taking as basis elements of the Schwartz space 
$q^j e^{-\frac{q^2}{2}} \in \mathcal{S}(\mR)$, $j\in \mN_0$, 
\begin{eqnarray}
e^{\alpha X_s}q^j e^{-\frac{q^2}{2}}=\sum_{k=0}^\infty 
\frac{\alpha^k}{k!} X_s^k q^j e^{-\frac{q^2}{2}}
\end{eqnarray} 
is a Schwartz function in $q$, because in the expansion of
$X_s^k q^j e^{-\frac{q^2}{2}}$ the maximal exponent of $q$ 
is just $k+j$, cf. \eqref{Xs^keminusqpul}.
Therefore, $e^{\alpha X_s}q^j e^{-\frac{q^2}{2}}$ grows as 
$q^je^{-\frac{q^2}{2}} e^{\alpha q}$, $\alpha\in \mC$, 
which is a characterizing property of Schwartz function class 
in the variable $q$.

\hfill

It is easy to verify the following identities in the universal enveloping 
algebra $\rm{U}\big(\sl(2, \mC)\big)$:
\begin{eqnarray}\label{DXnak}
\, [E+n,X_s^k] &=& k X_s^k,
\nonumber \\
\, [D_s,X_s^k] &=& -\mathrm{i}(E+n)X_s^{k-1}- \mathrm{i}X_s(E+n)X_s^{k-2}-\ldots  
- \mathrm{i}X_s^{k-1}(E+n) 
\nonumber \\
\, &=& -\mathrm{i} k\frac{k-1}{2}X_s^{k-1}- \mathrm{i} kX_s^{k-1}(E+n),
\end{eqnarray}
so that for all $\alpha\in \mC$
\begin{eqnarray}
&& [D_s,e^{\alpha X_s}]=\sum_{k=0}^\infty \frac{\alpha^k}{k!}[D_s,X_s^k] \nonumber \\
&& = -\mathrm{i}\sum_{k=1}^\infty \frac{\alpha^k}{(k-1)!}X_s^{k-1}(E+n)-
\mathrm{i}\frac{\alpha^2}{2}X_s\sum_{k=2}^\infty \frac{\alpha^{k-2}}{(k-2)!}X_s^{k-2} \nonumber\\
&& = -\mathrm{i}\alpha e^{\alpha X_s}(E+n)-\mathrm{i}\frac{\alpha^2}{2}X_s e^{\alpha X_s}.
\end{eqnarray}

The substitution of 
$$
f=e^{\alpha X_s} g(X_s)m_k^s
$$
into \eqref{eigensf}, where $m_k^s\in M_k^s$ is a symplectic monogenic and 
$g(X_s)$ is a polynomial in $X_s$, yields
$$
D_s e^{\alpha X_s} g(X_s)m_k^s-c X_se^{\alpha X_s} g(X_s)m_k^s
=\lambda e^{\alpha X_s} g(X_s)m_k^s.
$$
Because $e^{\alpha X_s}$ is an invertible operator, we get
$$
D_s (g(X_s)m_k^s)-\mathrm{i} \alpha (E+n) g(X_s)m_k^s - 
\big(c+\mathrm{i}\frac{\alpha^2}{2}\big)X_s g(X_s)m_k^s=\lambda g(X_s)m_k^s .
$$
Now we set $c=-\mathrm{i}\frac{\alpha^2}{2}$, i.e., $\sqrt{2\mathrm{i}c}=\alpha$ 
(we choose and fix one of the roots):
\begin{equation}\label{Dsg(Xs)}
D_s \big(g(X_s)m_k^s\big)-\mathrm{i}\alpha(E+n)g(X_s)m_k^s=\lambda g(X_s)m_k^s
\end{equation}
and substitute 
\begin{equation}
g(X_s)=g_k^j(X_s)=\sum_{l=0}^j \beta_l^{j,k}X_s^l.
\end{equation}
Then \eqref{Dsg(Xs)} turns into the recursion relation 
\begin{eqnarray*}
&& \lambda \sum_{l=0}^j \beta_l^{j,k}X_s^l m_k^s=
-\mathrm{i}\alpha\sum_{l=0}^j \beta_l^{j,k}(l+k+n)X_s^l m_k^s 
+\sum_{l=0}^j \beta_l^{j,k} D_s (X_s^l m_k^s),
\end{eqnarray*}
and noting $D_s (X_s^l m_k^s)=-\mathrm{i}\frac{l}{2}(2k+2n+l-1)X_s^{l-1} m_k^s$, see \eqref{DXnak}, we have
$$
\sum_{l=0}^j (\lambda +\mathrm{i}\alpha(l+k+n)) \beta_l^{j,k}X_s^l m_k^s = 
-\mathrm{i}\sum_{l=0}^{j-1}\frac{l+1}{2}(2k+2n+l)\beta_{l+1}^{j,k}X_s^l m_k^s. 
$$
Finally, we obtain the recurrence relations for $l=0,1,\ldots,j-1$:
\begin{eqnarray}
&& (\lambda +\mathrm{i}\alpha(l+k+n)) \beta_l^{j,k} 
= -\mathrm{i}\frac{l+1}{2}(2k+2n+l)\beta_{l+1}^{j,k},\\
&& (\lambda +\mathrm{i}\alpha(l+k+n)) \beta_j^{j,k} =0.
\end{eqnarray}
In order for $g(X_s)$ to be a polynomial in $x,y$ of degree $j$, we need to have 
$\lambda =-\mathrm{i}\alpha( n+j+k)$ as an eigenvalue. 
Hence our recursion becomes
$$
\alpha(j-l)\beta_l^{j,k} =\frac{l+1}{2} (2k+2n+l) \beta_{l+1}^{j,k}, 
$$
which results in
$$ 
\beta_{l+1}^{j,k} = \frac{2\alpha(j-l)}{(l+1)(2k+2n+l)}\beta_l^{j,k}=
\ldots 
=2^{l+1}\alpha^{l+1}\binom{j}{l+1}\frac{(2k+2n-1)!}{(2k+2n+l)!} \beta_0^{j,k}.
$$
Therefore, we conclude that
\begin{equation}
\beta_{l}^{j,k} = 2^{l}\alpha^l \binom{j}{l}\frac{(2k+2n-1)!}{(2k+2n-1+l)!} \beta_0^{j,k},
\end{equation}
and we choose $\beta_0^{j,k}=1$. Hence we have
\begin{eqnarray}
&& g_k^j(X_s)=\sum_{l=0}^j 2^l \alpha^l \binom{j}{l}\frac{(2k+2n-1)!}{(2k+2n-1+l)!}X_s^l \nonumber \\
&& =j! \frac{(2k+2n-1)!}{(2k+2n-1+j)!}L^{2n+2k-1}_j(-2\alpha X_s),
\end{eqnarray}
where $L_j^\beta$ is the generalized Laguerre polynomial,
$$
L_j^\beta (x)=\sum_{l=0}^j(-1)^l \frac{(j+\beta)(j+\beta-1)\ldots(j+\beta +l-j+1)}{(j-l)!}\frac{x^l}{l!}
$$
defined by the formula
$$
L_j^\beta(x)=\frac{x^{-\beta} e^{x}}{j!}\frac{d^j}{d x^j}(x^{x+\beta} e^{-x}).
$$
The spectral decomposition of our operator, which can be termed 
the symplectic spin harmonic oscillator, is summarized in the following theorem.
\begin{theorem}
The operator $H=D_s-cX_s$, $c\in \mR$, has a complete system of 
eigenfunctions (valued in the Segal-Shale-Weil representation) given by
\begin{equation}
f_k^j=e^{\sqrt{2\mathrm{i}c}X_s}L_j^{2n+2k-1}(-2\sqrt{2\mathrm{i}c}X_s ) m_k^s,
\end{equation}
where $L_j^\alpha(-2\sqrt{2\mathrm{i}c}X_s)$ is the generalized Laguerre polynomial of the operator 
$-2\sqrt{2\mathrm{i}c}X_s$, $m_k^s\in M_k^s$ is a symplectic monogenic and $j,k\in \mathbb{N}$, 
with corresponding eigenvalue
\begin{equation}
\lambda_k^j = \sqrt{2\mathrm{i}c}(n+j+k).
\end{equation}
\end{theorem}
\begin{example}
The simplest eigenfunction for $j=0$ 
is $e^{\sqrt{2\mathrm{i}c}X_s}e^{-\frac{q^2}{2}}\in M_0^s$, 
where $e^{-\frac{q^2}{2}}$ 
is a highest weight vector of the Segal-Shale-Weil representation.
\end{example}


\section{Fischer product and reproducing kernel on symplectic spinors}

Let us briefly mention a motivation given by the classical orthogonal 
Fischer scalar product. For two complex polynomials valued in the 
Clifford algebra associated to a quadratic form,
$f\otimes a$, $g\otimes b \in \mathrm{Pol}(\mR^m,\mC) \otimes Cl(\mR^m)$, the Fischer scalar product is defined by
$$\langle f\otimes a, g\otimes b\rangle = [\overline{f(\partial_x)}g]_{x=0}[\overline{a}b]_0.$$
Here $\overline{f(\partial_x)}$ is a differential operator, where we substitute $\partial_{x_j}$ 
for the variable $x_j$, $j=1,\ldots ,m$, and act by the resulting differential operator 
on a polynomial $g(x)$. As for the values, $[\,\,]_0$ denotes
the zero degree part of an element in $Cl(\mR^m)$.
The properties of scalar products are conveniently encoded in their reproducing kernels. For example, the space of homogeneous polynomials of homogeneity $k$ satisfies
$$
\left\langle \frac{\langle x,y \rangle^k}{k!} ,g(x)\right\rangle= g(y)
$$
for all $g \in \mathrm{Pol}_k(\mR^m,\mC)$ and $\langle , \rangle$ the canonical
scalar product on $\mR^m$. Hence the reproducing kernel for homogeneity $k$ 
harmonic polynomials $\mathcal{H}_k$,
\begin{equation}
Z_k(x,y)=\rm{Proj}_{\mathcal{H}_k} \left( \frac{\langle x,y \rangle^k}{k!}  \right),
\end{equation}
can be expressed by the use of the Gegenbauer polynomial.
The interested reader can find more about this topic in, e.g.,  \cite{DX} and \cite{dbws}.

In what follows, we attempt to apply the concept of Fischer product and 
reproducing kernel to the space of symplectic spinors equipped with the 
action of the metaplectic Lie algebra. As in the previous section, after 
some general considerations we focus mostly on the real dimension $2$. 

\subsection{Fischer product and reproducing kernel for $n=1$}

We now aim to define the Fischer product on the space of symplectic spinors. 
We construct the symplectic Fischer product on 
$\mathrm{Pol}(\mR^{2n},\mC)\otimes \mathcal{S}(\mR^n)$ for $f\otimes \psi$, $g\otimes \phi$ with $f,g\in \mathrm{Pol}(\mR^{2n},\mC)$ and $\psi, \phi \in \mathcal{S}(\mR^n)$, in the form
\begin{equation}\label{fisherproduct}
\langle f\otimes \psi , g\otimes \phi \rangle =\omega(f,g) \int_{\mR^n}\overline{\psi(q)}\phi(q) \,\mathrm{d}q.
\end{equation}
The integral is the inner product in the fiber variables $q_1, \ldots, q_n$ 
and $\omega(f,g)$ is the evaluation of a lift of the symplectic form to 
symmetric tensors $Sym_k(\mR^{2n})$, $k\in\mN$. We put
\begin{equation}\label{symplinnerprod}
\omega(v_1\otimes \ldots \otimes v_k,w_1\otimes \ldots \otimes w_k )
=\sum_{(j_1,\ldots, j_k)\in S_k} \omega(v_1,w_{j_1})\omega(v_2,w_{j_2})\ldots\omega(v_k,w_{j_k}),
\end{equation}
where $v_j,w_j\in \mR^{2n}$ and we sum over all even permutations of the set $\{1, \dots ,k\}$.

As already advertised above, we now focus on the real $2$-dimensional case and 
for a moment elaborate on the part of the inner product on $\mathrm{Pol}(\mR^{2},\mC)$
given in \eqref{symplinnerprod}.
We normalize the lift of the symplectic form to be $\omega(e_1,e_2)=1$ for 
$v=xe_1+ye_2\in\mR^2$, and define 
the Fourier symplectic transformation by 
$$
x\longleftrightarrow \partial_y, \quad y\longleftrightarrow - \partial_x. 
$$
Consequently, we get for $r,s,t,u\in\mN_0$
$$
\langle x^r y^s ,x^t y^u \rangle=\omega( x^r y^s ,x^t y^u )=
(-1)^s \partial_y^r \partial_x^s x^t y^u = (-1)^s u! s! \delta_{r,u} \delta_{s,t},
$$
and so we have for $f=x^r y^s, g=x^t y^u$ and $r+s=t+u$
\begin{eqnarray*}
&& \omega( f, y g)= \langle f, y g\rangle = 
(-1)^s \partial_y^r \partial_x^s x^t y^{u+1}=(-1)^s (u+1)! s! \delta_{r,u+1}\delta_{s,t}, \\
&& \omega( \partial_x f,g )=\langle \partial_x f,g \rangle = 
r \langle x^{r-1}y^s,x^t y^u \rangle =(-1)^s r (r-1)! s! \delta_{r-1,u}\delta_{s,t}, \\
&& \omega( f, x g)= \langle f, x g\rangle = (-1)^s (t+1)! s! \delta_{r,u}\delta_{s,t+1}, \\
&& \omega( \partial_y f,g )=\langle \partial_y f,g \rangle = 
(-1)^{s-1} r ! (s-1)! \delta_{r,u}\delta_{s-1,t}. 
\end{eqnarray*}
Hence, there are the relations
\begin{equation}
\langle \partial_x f,g \rangle = \langle f, yg \rangle,\quad  - \langle \partial_y f,g \rangle = \langle f, xg \rangle.
\end{equation}
Let us now summarize our definitions and basic properties in the $2$-dimensional case.
\begin{definition}
The symplectic Fischer product for $f(x,y)\otimes \psi$, $g(x,y)\otimes \phi$, with 
$f,g\in \mathrm{Pol}(\mR^2,\mC)$ and $\psi, \phi\in \mathcal{S}(\mR)$, is given by
\begin{equation}\label{sfsp}
\langle f\otimes \psi, g\otimes \phi \rangle = 
\left[ f(\partial_y,-\partial_x) g(x,y) \right]_{x=y=0} \int_{-\infty}^{\infty} \overline{\psi(q)} \phi(q) \,\mathrm{d}q,
\end{equation}
where the bar denotes the complex conjugation of a complex valued function.
\end{definition}

\begin{lemma}
The bilinear form defined in \eqref{sfsp} for all $a,b\in \mathrm{Pol}(\mR^2,\mC)\otimes {\mathcal{S}}(\mR)$ satisfies 
\begin{enumerate}
\item $\langle q a,b \rangle = \langle a, q b\rangle $, and 
$\langle \mathrm{i} q a,b \rangle = \langle a,-\mathrm{i}q b\rangle $.
\item $\langle \partial_q a,b \rangle = -\langle a, \partial_q b\rangle $ 
and $\langle \mathrm{i} \partial_q a,b \rangle = \langle a, \mathrm{i} \partial_q b\rangle $.
\item $\langle \partial_x a,b \rangle = \langle a, y b\rangle $.
\item $\langle \partial_y a,b \rangle = -\langle a, x b\rangle $.
\item $\langle x a,b \rangle = \langle a, \partial_y b\rangle $.
\item $\langle y a,b \rangle = -\langle a, \partial_x b\rangle $.
\end{enumerate}
\end{lemma}
Now we compute the adjoints of operators $D_s, X_s$ with respect to $\langle\, , \rangle $.
\begin{lemma}
The adjoint operator for the symplectic Dirac operator $D_s$ with respect to the symplectic Fischer 
product is $X_s$, and vice versa. We have
$$\langle D_s a,b \rangle=\langle  a, X_s b \rangle, \quad \langle X_s a,b \rangle=\langle a, D_s b \rangle,$$
for arbitrary $a,b\in \mathrm{Pol}(\mR^2,\mC)\otimes \mathcal{S}(\mR)$.
\end{lemma}

{\bf Proof:} A direct computation for $a,b\in \mathrm{Pol}(\mR^2,\mC)\otimes \mathcal{S}(\mR)$ gives
\begin{eqnarray*}
&& \langle D_s a,b \rangle = \langle (\mathrm{i} q \partial_y - \partial_q \partial_x)a,b\rangle= 
\langle a, (\mathrm{i} q x + \partial_q y)b \rangle,\\
&& \langle X_s a,b \rangle = \langle (\mathrm{i} q x + \partial_q y)a,b\rangle= 
\langle a, (\mathrm{i} q \partial_y - \partial_q \partial_x)b \rangle.
\end{eqnarray*}

\hfill
$\square$

Consequently, we have the orthogonality relations 
for the symplectic Fischer decomposition,
\begin{equation}
\langle X_s^j m_k^s , X_s^l m_h^s \rangle \sim \delta_{j,l} \delta_{k,h},
\end{equation}
with symplectic monogenics $m_k^s\in M_k^s, m_h^s\in M_h^s$.

\begin{lemma}\label{HXYadjoint}
The adjoint operators to the basis elements $ \tilde{X},\tilde{Y}$
and $\tilde{H}$ of $\mp(2,\mR)$, cf. \eqref{HXYxy}, with respect to 
the symplectic Fischer product are $ -\tilde{X},-\tilde{Y}$ and $-\tilde{H}$, 
respectively.
\end{lemma}
Now we pass to the construction of the reproducing kernel $K_k(\xi_1,\xi_2, x, y)$ 
for the bilinear form 
$$
(f(x,y),g(x,y))=\left[ f(\partial_y,-\partial_x) g (x,y)\right]_{x=y=0}
$$
on the space of polynomials of homogeneity $k$. 
Inspired by the orthogonal case, we claim 
\begin{equation}
K_k(\xi_1,\xi_2,x,y)=\frac{1}{k!}\big(-\xi_1 y + \xi_2 x \big)^k .
\end{equation}
Indeed, we have
\begin{eqnarray}\label{Kxixi}
&& \big(K_k(\xi_1,\xi_2,x,y),p(x,y)\big)=\left( \frac{1}{k!}\big(-\xi_1 y + \xi_2 x \big)^k,p(x,y) \right)\nonumber \\
&& =\frac{1}{k!}\big(\xi_1 \partial_x + \xi_2 \partial_y \big)^k p(x,y)=p(\xi_1,\xi_2)
\end{eqnarray}
for $p(x,y)\in \mathrm{Pol}_k(\mR^2,\mC)$.

In order to adapt $K_k(\xi_1,\xi_2,x,y)$ 
to the reproducing kernel $Z_k$ of the space of symplectic monogenics $M_k^s$,
we shall regard $K_k(\xi_1,\xi_2,x,y)$ as an element in the algebra
${\mathrm{Pol}_k(\mR^2\times\mR^2,\mC)\otimes\mathrm{End}(\mathcal{S}(\mR)})$
with the value in $\mathrm{End}(\mathcal{S}(\mR))$ 
given by the identity endomorphism on $\mathcal{S}(\mR)$.
Moreover, we introduce the projector
\begin{eqnarray}
& & \mathrm{Proj}_{sm}^k:\, \mathrm{Pol}_k(\mR^2,\mC)\otimes \mathcal{S}(\mR)\rightarrow M_k^s,
\nonumber \\
& & \mathrm{Proj}_{sm}^k = \sum_{j=1}^k a_j^k X_s^j D_s^j\,\,\in
{\mathrm{Pol}_k(\mR^2,\mC)\otimes \mathrm{End}(\mathcal{S}(\mR))},
\end{eqnarray}
to homogeneity $k$ symplectic monogenics, see \cite{bss}, and 
define the symplectic Fischer $\mathrm{End}(\mathcal{S}(\mR))$-valued 
pairing for the elements in the spaces
${\mathrm{Pol}_k(\mR^2\times\mR^2,\mC)\otimes \mathrm{End}(\mathcal{S}(\mR))}$
 and ${\mathrm{Pol}_k(\mR^2,\mC)\otimes \mathrm{End}(\mathcal{S}(\mR))}$ by
\begin{eqnarray}\label{fischer2}
\langle f(\xi_1,\xi_2,x,y,q,\partial_q),g(x,y,q)\rangle =
\left[ {f(\xi_1,\xi_2,\partial_y,-\partial_x,q,\partial_q)} g(x,y,q)\right]_{x=y=0}.
\end{eqnarray}
We remark that we used in \eqref{fischer2} the same notation $\langle\, , \rangle$ for
the symplectic Fischer $\mathrm{End}(\mathcal{S}(\mR))$-valued product
as for the $\mR$-valued scalar product \eqref{fisherproduct}, and believe 
the attentive reader will not have a problem in distinguishing
which of them is currently used. Another remark is that we exploit 
in \eqref{fischer2}
 the well-known fact that any symplectic spinor 
$g(x,y,q)\in {\mathrm{Pol}_k(\mR^2,\mC)\otimes \mathcal{S}(\mR)}$ can be
regarded as an element in ${\mathrm{Pol}_k(\mR^2,\mC)\otimes \mathrm{End}(\mathcal{S}(\mR))}$,
because the space of Schwartz functions is a complex algebra.

\begin{theorem} 
The projection operator $\mathrm{Proj}_{sm}^k$ and the reproducing kernel $Z_k$
relate to the symplectic Fischer product as follows:
\begin{enumerate}
\item $\mathrm{Proj}_{sm}^k$ is self-adjoint.
\item $Z_k(\xi_1,\xi_2,x,y,q,\partial_q)=\mathrm{Proj}_{sm}^k K_k(\xi_1,\xi_2,x,y)$ is the reproducing kernel for $M_k^s$.
\end{enumerate}
\end{theorem}
{\bf Proof:}
Indeed, using this pairing, we first observe the self-adjointness property of $\mathrm{Proj}_{sm}^k$: 
$$ \langle \mathrm{Proj}_{sm}^k f,g \rangle = \sum_{j=0}^k a_j^k \langle X_s^j D_s^j f,g \rangle 
 =\sum_{j=0}^ka_j^k \langle f, X_s^j D_s^j  g \rangle = \langle f, \mathrm{Proj}_{sm}^k g \rangle . $$
By \eqref{Kxixi}, we have for $m_k^s\in M_k^s$
\begin{eqnarray*}
&& \langle Z_k(\xi_1,\xi_2,x,y,q,\partial_q), m_k^s(x,y,q) \rangle = \langle K_k(\xi_1,\xi_2,x,y), \mathrm{Proj}_{sm}^k m_k^s \rangle \\
&& = m_k^s(\xi_1,\xi_2,q),
\end{eqnarray*}
and for any $j\in\mN$ holds
\begin{eqnarray*}
\langle Z_k(\xi_1,\xi_2,x,y,q,\partial_q),X_s^j m_{k-j}^s(x,y,q) \rangle
= \langle K_k(\xi_1,\xi_2,x,y), \mathrm{Proj}_{sm}^k X_s^j m_{k-j}^s \rangle =0 .
\end{eqnarray*}

\hfill
$\square$

\begin{proposition}
The reproducing kernel $Z_k$ has the explicit form 
\begin{eqnarray}
Z_k(\xi_1,\xi_2,x,y,q,\partial_q)=\sum_{j=0}^k {\mathrm{i}}^j a_j^k \frac{1}{(k-j)!}
\big(-\xi_1 y  + \xi_2 x \big)^{k-j} X_s^j \xi_s^j, 
\end{eqnarray}
where $\xi_s=-q \xi_1 + \mathrm{i} \partial_q \xi_2$.
\end{proposition}
{\bf Proof:}
First we need an explicit formula for $D_s^j K_k(\xi_1,\xi_2,x,y)$, $j=1,\ldots,k$. 
We obtain by the chain rule
\begin{eqnarray*}
&& D_s K_k(\xi_1,\xi_2,x,y)=(\mathrm{i} q \partial_y - \partial_q \partial_x)\frac{1}{k!}
\big(-\xi_1 y + \xi_2 x \big)^k\\
&&=  K_{k-1}(\xi_1,\xi_2,x,y)(\mathrm{i}q \xi_1-\partial_q \xi_2)=
\mathrm{i} K_{k-1}(\xi_1,\xi_2,x,y)\xi_s.
\end{eqnarray*}
Therefore, $D_s^j K_k(\xi_1,\xi_2,x,y)=\mathrm{i}^j K_{k-j}(\xi_1,\xi_2,x,y) \xi_s^j$, and so
\begin{eqnarray*}
&& Z_k(\xi_1,\xi_2,x,y,q,\partial_q)= \sum_{j=0}^k a_j^k X_s^j D_s^j K_k(\xi_1,\xi_2,x,y)\\
&& = \sum_{j=0}^k \mathrm{i}^j a_j^k X_s^j K_{k-j}(\xi_1,\xi_2,x,y) \xi_s^j=
\sum_{j=0}^k \mathrm{i}^j a_j^k X_s^j \frac{\big(-\xi_1 y  + \xi_2 x \big)^{k-j}}{(k-j)!}  \xi_s^j
\end{eqnarray*}
which proves the assertion.

\hfill
$\square$

\section{Explicit bases of symplectic monogenics on $(\mR^{2},\omega)$}

In the present section we construct some explicit bases for symplectic monogenics 
$M_h^s$ in $\mathrm{Pol}(\mR^2,\mC)\otimes \mathcal{S}(\mR)$ of homogeneity $h$, and 
prove several useful characterizing properties.

The first distinguished basis is written in the real coordinates $x$ and $y$ on
$\mR^2$ and the (topological) basis $q^j e^{-\frac{q^2}{2}}$, $j\in \mN_0$,
of the Schwartz space. The second distinguished basis for symplectic monogenics 
is written in the complex coordinates $z$ and $\bar{z}$ on $\mR^2\simeq\mC$ 
and the (topological) basis of Hermite functions $\psi_j(q)$, $j\in \mN_0$, 
for $\mathcal{S}(\mR)$.
\begin{proposition}\label{basissmhol}
The symplectic spinors of homogeneity $h$ in the variables $x,y$ and odd in the 
variable $q$ for $h,k\in \mN_0$ with $k\geq h$,
\begin{equation}\label{basexyo}
\tilde{s}_{o,k}^h = e^{-\frac{q^2}{2}}\sum_{p=0}^h 
(-1)^p \frac{(2k+1)!! }{(2k-2p+1)!!}\binom{h}{p} q^{2k+1-2p}(x+\mathrm{i} y)^{h-p}(\mathrm{i}y)^p
\end{equation}
and even in variable $q$ for $k\in \mN_0$,
\begin{equation}\label{basexye}
\tilde{s}_{e,k}^h = e^{-\frac{q^2}{2}}\sum_{p=0}^h (-1)^p \frac{(2k)!! }{(2k-2p)!!}
\binom{h}{p} q^{2k-2p}(x+\mathrm{i} y)^{h-p}(\mathrm{i} y)^p,
\end{equation}
form a (topological) basis of the odd and even part of the symplectic monogenics $M_h^s$, respectively.
\end{proposition}

{\bf Proof:}
Let us consider a polynomial monogenic symplectic spinor 
\begin{equation*}
f(x,y,q)=e^{-\frac{q^2}{2}}  \sum_{j=0}^{\infty} q^j p_j(x,y),
\end{equation*}
where $p_j(x,y)$ are polynomials in the variables $x, y$. 
Solving the equation $D_s f(x,y,q)=0$, we have
\begin{eqnarray*}
0 & =&(\mathrm{i}q\partial_y-\partial_x \partial_q) f(x,y,q)\\
& =&e^{-\frac{q^2}{2}}  \sum_{j=0}^{\infty} 
\left( \mathrm{i}q^{j+1}\partial_y p_j(x,y)+q^{j+1} 
\partial_x p_j(x,y)-j q^{j-1}\partial_x p_j(x,y)\right).
\end{eqnarray*}
The Schwartz functions $e^{-\frac{q^2}{2}}q^j$, $j\in \mN_0$ are linearly independent, hence 
\begin{equation}
q^j\left( (\partial_x+i\partial_y ) p_{j-1}(x,y)-(j+1)\partial_x p_{j+1}(x,y) \right)=0
\end{equation}
for each $j\in \mN_0$. We get a system of recursion equations, splitting into two subsystems 
of odd and even homogeneity in the variable $q$ and the solution follows.

For a fixed homogeneity $h$ in the variables $x$ and $y$, the systems $\tilde{s}_{o,k}^h$ 
and $\tilde{s}_{e,k}^h $ contain all powers $q^j$, $j\in \mN_0$, for appropriate $k$ and 
all possible combinations of $x,y$ in $\mathrm{Pol}(\mR^{2},\mC)$ so that they are in $\Ker(D_s)$. 
Therefore, the odd \eqref{basexyo} and even \eqref{basexye} systems form a basis of 
$M_h^s$ because $\{q^j e^{-\frac{q^2}{2}}\}_{j\in\mN}$ is a (topological) 
basis of $\mathcal{S}(\mR)$.

\hfill
$\square$

In the complex coordinates $z=x+iy,\overline{z}=x-iy$,
with $\partial_x=\partial_z+\partial_{\bar{z}}$ and 
$\partial_y=i(\partial_z-\partial_{\bar{z}})$, the 
symplectic Dirac operator is by \eqref{XsHolom} given by 
\begin{equation}\label{DsHolom2} 
D_s=-(q+\partial_q)\partial_z+(q -\partial_q)\partial_{\bar{z}}.
\end{equation}
Let us recall (see e.g. \cite{Sz}) that the Hermite functions $\{\psi_k(q)\}_{k\in\mN_0}$ 
form a (topological) basis of the Schwartz space $\mathcal{S}(\mR)$.
The $k$-th Hermite function is
$$\psi_k(q)=\frac{1}{\sqrt{2^k k! \sqrt{\pi}}}e^{-\frac{q^2}{2}} H_k(q)=\frac{(-1)^k}{\sqrt{2^k k! \sqrt{\pi}}}\big(q-\partial_q\big)^k e^{-\frac{q^2}{2}},$$
where $H_k$ is the $k$-th Hermite polynomial.
The operators $(q+\partial_q)$ and $(q-\partial_q)$ act on the basis vectors by
\begin{eqnarray}\label{vzorceqDq}
(q+\partial_q) \psi_k&=&\sqrt{2}\sqrt{k}\psi_{k-1},\nonumber\\
(q-\partial_q)\psi_k&=& \sqrt{2}\sqrt{k+1} \psi_{k+1},
\end{eqnarray}
and together with the operator acting by a multiple of identity on each 
$\psi_k$ form the representation of the Lie algebra $\sl(2,\mC)$. We shall use   
the following easily verified formulas
\begin{eqnarray}\label{vzorceqDq2}
(q^2-\partial_q^2)\psi_k &=& (2k+1)\psi_k,\nonumber\\
(q-\partial_q)^2 \psi_k &=& 2\sqrt{(k+1)(k+2)}\psi_{k+2},\nonumber\\
(q+\partial_q)^2 \psi_k &=& 2\sqrt{k(k-1)}\psi_{k-2}.
\end{eqnarray}
\begin{proposition}\label{lemHerm}
The polynomial symplectic spinors of homogeneity $h$ in the variables $z, \bar{z}$ 
and odd in the variable $q$ for $k\in  \mathbb{N}_0$,
\begin{equation}\label{hermOdd}
s_{o,k}^h=\sum_{p=0}^h \sqrt{\frac{(2k+2p)!!}{(2k+2p+1)!!}}\binom{h}{p} \psi_{2k+2p+1}(q)\bar{z}^{h-p} z^p,
\end{equation}
form the basis of the odd part of the solution space of the symplectic Dirac operator $D_s$.

The polynomial symplectic spinors of homogeneity $h$ in the variables $z, \bar{z}$ and even in the variable $q$ for 
$k\in \mathbb{N}_0$,
\begin{equation}\label{hermEvenK}
s_{e,k}^h=\sum_{p=0}^h \sqrt{\frac{(2k+2p-1)!!}{(2k+2p)!!}}\binom{h}{p} \psi_{2k+2p}(q)\bar{z}^{h-p} z^p,
\end{equation}
and for $k=-1, -2,\ldots, -h$
\begin{equation}\label{hermEven0}
s_{e,k}^h=\sum_{p=|k|}^h \sqrt{\frac{(2k+2p-1)!!}{(2k+2p)!!}}\binom{h}{p} \psi_{2k+2p}(q)\bar{z}^{h-p} z^p,
\end{equation}
form the basis of the even part of the solution space of the symplectic Dirac operator $D_s$.
\end{proposition}
{\bf Proof:}
Let us consider a polynomial monogenic symplectic spinor 
$$f(z,\bar{z},q)=\sum_{l=0}^{\infty} \psi_l(q) p_l(z,\bar{z}),$$
where $\psi_l(q)$ is the $l$-th Hermite function and $p_l(z,\bar{z})$ is a polynomial in the variables $z,\bar{z}$.
The action of the symplectic Dirac operator is then  
\begin{eqnarray*}
0=D_s f(z,\bar{z},q)&=& \big((q+\partial_q)\partial_z-(q -\partial_q)\partial_{\bar{z}} \big)f(z,\bar{z},q)\\
&=&\sqrt{2}\sum_{l=0}^\infty \sqrt{l}\psi_{l-1}(q)\partial_z p_l (z,\bar{z})-\sqrt{l+1}\psi_{l+1}(q)\partial_{\bar{z}}p_l(z,\bar{z}).
\end{eqnarray*}
The Hermite functions are linearly independent, which implies
\begin{equation}
\psi_l(q)(\sqrt{l+1}\partial_z p_{l+1}(z,\bar{z})-\sqrt{l}\partial_{\bar{z}}p_{l-1}(z,\bar{z}))=0
\end{equation}
for each $l\in \mathbb{N}_0$. The system of recursion equations is split into two systems 
with odd and even indexes in the variable $q$, each of which is easy to resolve.

For a fixed homogeneity $h$ the systems of symplectic polynomial spinors \eqref{hermOdd}, \eqref{hermEvenK} and \eqref{hermEven0} form a basis of symplectic monogenics $M_k^s$ of homogeneity $h$, 
and because the Hermite functions form a basis of $\mathcal{S}(\mR)$ the above collection of 
symplectic monogenics is a (topological) basis of $\Ker(D_s)$.

\hfill
$\square$

Let us now explore the properties of the symplectic Fischer product \eqref{sfsp} applied to the
basis elements discussed in the Proposition \ref{basissmhol} and Proposition \ref{lemHerm}.
The motivation for this question is the existence of a basis of symplectic monogenics, 
which is isotropic with respect to the product \eqref{sfsp}.

\begin{lemma}
The basis elements \eqref{basexyo} and \eqref{basexye} 
of homogeneity $2$ in the symplectic Fischer product \eqref{sfsp} 
satisfy, for $k,l\in \mN$, $k,l\geq 2$,
\begin{eqnarray*}
&& \langle \tilde{s}_{o,k}^2,\tilde{s}_{o,l}^2 \rangle = \frac{-3\sqrt{\pi}(2k+2l-5)}{2^{k+l-3}},\\
&& \langle \tilde{s}_{e,k}^2,\tilde{s}_{e,l}^2 \rangle = \frac{-3\sqrt{\pi}(2k+2l-5)}{2^{k+l-2}},\\
&& \langle \tilde{s}_{o,k}^2,\tilde{s}_{e,l}^2 \rangle =0.
\end{eqnarray*}
\end{lemma}
{\bf Proof:}
Focusing just on the $\mathrm{Pol}(\mR^2,\mC)$ part of the symplectic Fischer product  
\eqref{sfsp}, the only non zero combinations of $x,y$ in homogeneity 2 are 
$\langle x^2,y^2 \rangle=2$ and $\langle xy,xy \rangle=-1$.
Then $\int_{-\infty}^{\infty} e^{-q^2}q^{2t}\, \,\mathrm{d}q=\frac{\sqrt{\pi}(2t+1)}{2^t}$ 
for $t\in \mN_0$ and moreover, $\int_{-\infty}^{\infty} e^{-q^2}q^{t}\,\mathrm{d}q=0$ 
for $t$ odd.

\hfill
$\square$

Therefore, we see that the symplectic Fischer product \eqref{sfsp} of any 
two odd or even basis elements \eqref{basexyo}, \eqref{basexye} 
for $k,l\geq 2$ is non-zero (in fact, negative) for $k=l$. 
This implies that the symplectic Fischer product \eqref{sfsp} does not
seem to be a convenient candidate for the scalar product on 
$\mathrm{Pol}(\mR^2,\mC)\otimes \mathcal{S}(\mR)$. 

Let us rewrite the symplectic Fischer product \eqref{sfsp} in the complex 
variables. In the variables $z, \bar{z}$, we have a non-trivial pairing for
the pairs
$z\longleftrightarrow -2\mathrm{i}\partial_{\bar{z}}$ and $\bar{z} 
\longleftrightarrow 2\mathrm{i} \partial_z$. 
Hence for $f(z,\bar{z})\otimes \psi$, $g(z,\bar{z})\otimes \phi$, with 
$f,g\in \mathrm{Pol}(\mR^2,\mC)$ and $\psi, \phi\in \mathcal{S}(\mR)$,
\begin{equation}\label{sfspzz}
\langle f\otimes \psi, g\otimes \phi \rangle = 
\left[ f(-2\mathrm{i}\partial_{\bar{z}},2\mathrm{i}\partial_z) 
g(z,\bar{z}) \right]_{z=\bar{z}=0} \int_{-\infty}^{\infty} 
\overline{\psi(q)} \phi(q) \,\mathrm{d}q .
\end{equation}
Let us look at the symplectic Fischer product \eqref{sfspzz} for the 
low homogeneity basis elements $s_{o,k}^h$, \eqref{hermOdd}, of odd 
part of the symplectic monogenics. 
\begin{example}
In the homogeneity $2$ and $k,l\in \mN$ holds
\begin{eqnarray*}
\langle s_{o,k}^2, s_{o,l}^2\rangle &=& 
-8 \frac{(2k)!!}{(2k+1)!!}\delta_{2k+1,2l+5} +16\frac{(2k+2)!!}{(2k+3)!!} \delta_{2k+2,2l+2}\\
&& -8 \frac{(2k+4)!!}{(2k+5)!!}\delta_{2k+5,2l+1},
\end{eqnarray*}
where just one of the Kronecker deltas on the previous display 
may be non-zero. We observe that for $k=l$ holds
$\langle s_{o,k}^2, s_{o,k}^2\rangle\neq 0$, because $\delta_{2k+2,2l+2}\not=0$.

In the homogeneity $3$ and $k,l\in \mN$,
\begin{eqnarray*}
 \langle s_{o,k}^3,s_{o,l}^3 \rangle &=& -48 i \frac{(2k)!!}{(2k+1)!!}\delta_{2k+1,2l+7} -16 i \frac{(2k+2)!!}{(2k+3)!!}\delta_{2k+3,2l+5}\\
&& +16i \frac{(2k+4)!!}{(2k+5)!!}\delta_{2k+5,2l+3}+48i \frac{(2k+6)!!}{(2k+7)!!} \delta_{2k+7,2l+1},
\end{eqnarray*}
where again just one Kronecker delta may be non-zero. 
For $k=l$ the symplectic Fischer product gives zero, $\langle s_{o,k}^3,s_{o,k}^3 \rangle=0$,
and the analogous conclusion $\langle s_{o,k}^h, s_{o,k}^h\rangle=0$ can be made for 
all odd homogeneities $h$.
\end{example}

Let us now consider another skew-symmetric bilinear form on 
$\mathrm{Pol}(\mR^2,\mC)\otimes \mathcal{S}(\mR)$, which is skew-symmetric 
on $\mathcal{S}(\mR)$ and possesses several remarkable properties. We use again the
complex variables on ${\mathbb R}^2$. 
\begin{definition}
Let us introduce a bilinear form $\langle\, ,\rangle_1$ on symplectic spinors, 
defined on $f(z,\bar{z})\otimes \psi$, $g(z, \bar{z})\otimes \phi$ with 
$f,g\in \mathrm{Pol}(\mR^2,\mC)$ and $\psi, \phi\in \mathcal{S}(\mR)$ by
\begin{equation}\label{sympFormDef}
\langle f\otimes \psi, g\otimes \phi \rangle_1 = \sqrt{2}
\left[\frac{1}{h!} f(\partial_z,\partial_{\bar{z}}) g(z,\bar{z}) \right]_{z=\bar{z}=0} 
\int_{-\infty}^{\infty} \big( \partial_q\psi(q)\big) \phi(q) \,\mathrm{d}q,
\end{equation}
where $h$ denotes the homogeneity of the polynomial $f(z,\bar{z})$.
\end{definition}
In the monomial basis, we have for $r,s,t,u\in\mN_0$
\begin{equation}\label{sympFromZbqrZ}
\langle z^r \bar{z}^s\otimes \psi ,z^t \bar{z}^u\otimes \phi \rangle_1=
\sqrt{2}\frac{r! s!}{(r+s)!}\delta_{r,t} \delta_{s,u} \int_{-\infty}^{\infty} 
\big( \partial_q\psi(q)\big) \phi(q) \,\mathrm{d}q,
\end{equation}
where $\delta_{r,t}$ denotes the Kronecker delta. Moreover, for 
$a,b\in \mathrm{Pol}(\mR^2,\mC)\otimes \mathcal{S}(\mR)$ holds
\begin{eqnarray}
 \langle z a,b \rangle_1 = \langle a,\partial_z b \rangle_1,   
 && \langle \partial_{z} a,b \rangle_1 = \langle a, z b \rangle_1,\nonumber  \\
 \langle \bar{z} a,b \rangle_1 = \langle a,\partial_{\bar{z}} b \rangle_1, 
&& \langle \partial_{\bar{z}} a,b \rangle_1 = \langle a, \bar{z} b \rangle_1.
\end{eqnarray}
Notice that the bilinear form $\langle\, ,\rangle_1$ is not $\mp(2,\mR)$-invariant
on the whole space of symplectic spinors, because
\begin{eqnarray*}
&&\langle H_t (f\otimes \psi), g\otimes \phi \rangle_1 
- \langle f\otimes \psi, H_t (g\otimes \phi) \rangle_1\\
&&\quad\quad = [f,g]\int_{-\infty}^{\infty} q \psi(q) \phi(q) \,\mathrm{d}q,\\
&&\langle X_t (f\otimes \psi), g\otimes \phi \rangle_1 
- \langle f\otimes \psi, X_t (g\otimes \phi) \rangle_1\\
&&\quad\quad = \frac{\mathrm{i}}{2}[f,g]\int_{-\infty}^{\infty} \Big(q \psi(q) \phi(q) 
+2q \big(\partial_q  \psi(q)\big)  \big(\partial_q \phi (q)\big) \Big) \,\mathrm{d}q,\\
&&\langle Y_t (f\otimes \psi), g\otimes \phi \rangle_1 - \langle f\otimes \psi, Y_t (g\otimes \phi) \rangle_1\\
&&\quad\quad = \frac{\mathrm{i}}{2}[f,g]\int_{-\infty}^{\infty} \Big(q \psi(q) \phi(q) 
-2q \big(\partial_q  \psi(q)\big)  \big(\partial_q \phi (q)\big) \Big) \,\mathrm{d}q,
\end{eqnarray*}
with
\begin{eqnarray*}
 [f,g]:=\sqrt{2}\left[\frac{1}{h!} f(\partial_z,\partial_{\bar{z}}) g(z,\bar{z}) \right]_{z=\bar{z}=0}.
\end{eqnarray*}
However, $\langle\, ,\rangle_1$ is $\mp(2,\mR)$-invariant when restricted to 
any of the two irreducible subspaces of symplectic spinors (given by the subspaces 
of even and odd Schwartz functions, respectively.)

Let us now define the elements
\begin{equation}
s_{E,l}^h=\frac{1}{2^h}\sum_{j=0}^l s_{e,j}^h,
\end{equation}
which form a basis of even symplectic monogenics 
in the homogeneity $h$ (as well as the set 
$s_{e,l}^h$, $l\in \mN_0$, cf. \eqref{hermEvenK}.)

\begin{theorem}
The basis elements $s_{o,k}^h$, $s_{E,k}^h$, $k\in \mN_0$, of the polynomial 
symplectic monogenics of homogeneity $h$ in the variables $z,\bar{z}$ form 
two isotropic subspaces of the symplectic monogenics $M_h^s$ with respect to 
the form defined in \eqref{sympFormDef}. Namely, the basis elements satisfy
\begin{eqnarray}\label{isotropybase}
 \langle s_{o,k}^h, s_{o,l}^h  \rangle_1 =0, 
&& \langle s_{o,k}^h, s_{E,l}^h  \rangle_1 = \delta_{k,l}, \nonumber\\
 \langle s_{E,k}^h, s_{E,l}^h  \rangle_1 =0,
&&   \langle  s_{E,l}^h,s_{o,k}^h  \rangle_1 = -\delta_{k,l},
\end{eqnarray}
for $k,l\in \mN_0$ and $h\in \mN_0$. Moreover, the form is identically zero 
for symplectic monogenics of different homogeneities $h,h'$. 
\end{theorem} 
{\bf Proof:}
Let us remind the orthonormality relation
$\int_{-\infty}^{\infty} \psi_k (q) \psi_l (q)\,\mathrm{d}q = \delta_{k,l}$ 
for Hermite functions. Then the relations in the first column \eqref{isotropybase} are 
obvious, because the derivative of a Hermite function $\psi_k(q)$ is 
$$\partial_q \psi_k(q)=\sqrt{\frac{k}{2}}\psi_{k-1}(q)-\sqrt{\frac{k+1}{2}} \psi_{k+1}(q)$$
and consequently, the integral in the bilinear form is zero.

As for the proof of the relation 
$\langle s_{o,k}^h, s_{E,l}^h \rangle_1 = \delta_{k,l}$, we first prove 
$$\langle s_{o,k}^h, s_{e,l}^h \rangle_1 = 2^h (\delta_{k,l}- \delta_{k+1,l})$$
for $k,l\in \mN_0$. 
We use \eqref{sympFromZbqrZ} to simplify the calculation and get
\begin{eqnarray*}
 \langle s_{o,k}^h, s_{e,l}^h \rangle_1  &=& 
\sum_{p=0}^h \sqrt{\frac{(2k+2p)!!(2l+2p-1)!!}{(2k+2p+1)!!(2l+2p)!!}}\binom{h}{p}^2 
\frac{(h-p)!p!}{h!}\times \\
&& \big( \delta_{2k+2p,2l+2p} \sqrt{2k+2p+1} - \delta_{2k+2p+2, 2l+2p} \sqrt{2k+2p+2}  \big).
\end{eqnarray*}
This is equal to $\sum_{p=0}^h \binom{h}{p}=2^h$ for $k=l$, $-2^h$ 
for $k+1=l$ and zero otherwise. Then for the basis elements $s_{E,l}^h$ we have 
$\langle s_{o,k}^h, s_{E,l}^h \rangle_1= \sum_{j=0}^l \delta_{k,j}- \sum_{j=0}^l \delta_{k+1,j}$, 
which is non-zero just for $k=l$. 
The last relation in \eqref{isotropybase} follows from the skew-symmetry of the integration
in the variable $q$.
For different homogeneities, the statement easily follows from \eqref{sympFromZbqrZ}.

\hfill
$\square$

We remark that for $k<0$, the elements $s_{e,k}^h$ in \eqref{hermEven0} satisfy
$$\langle s_{e,l}^h,s_{e,k}^h\rangle_1=0, \quad  \langle s_{o,j}^h,s_{e,k}^h\rangle_1=0,\quad  \langle s_{E,j}^h,s_{e,k}^h\rangle_1=0,$$
for each  $k,l\in \mZ$, $-h \leq k< 0$, $-h\leq l$ and  $j\in \mN_0$. 

\subsection{The action of symmetry operators of the symplectic Dirac operator 
on the basis of symplectic monogenics}

In the present part we determine the action of the symmetry operators introduced in Section 
\ref{opervzbarz} on the basis of symplectic monogenics described in 
Proposition \ref{lemHerm}. We remark that the action of the symmetry operators
on the basis of symplectic monogenics described in Proposition \ref{basissmhol}
is much more involved.

We shall start with the even component of the basis.

\begin{proposition}
The operators $\partial_z$ and $\partial_{\bar{z}}$ decrease the homogeneity in  
$z,\bar{z}$ and preserve the elements of even basis \eqref{hermEvenK}, \eqref{hermEven0} 
of the kernel of the symplectic Dirac operator $D_s$. In particular, for $k\in \mZ$, $k\leq -h$,
\begin{eqnarray}
&&\partial_z s^h_{e,k} = h s^{h-1}_{e,k+1},\nonumber\\
&&\partial_{\bar{z}} s^h_{e,k} = h s^{h-1}_{e,k},\quad \mbox{for } k\neq -h, \quad \partial_{\bar{z}} s^h_{e,-h} = 0.
\end{eqnarray}
\end{proposition}
{\bf Proof}:
We verify the relation $\partial_z s^h_{e,k} = h s^{h-1}_{e,k+1}$ for $k\in \mN_0$, 
the others being analogous. We have
\begin{eqnarray*}
&&\partial_z\sum_{p=0}^h \sqrt{\frac{(2k+2p-1)!!}{(2k+2p)!!}}\binom{h}{p} \psi_{2k+2p}(q)\bar{z}^{h-p} z^p\\
&&=\sum_{p=1}^h \sqrt{\frac{(2k+2p-1)!!}{(2k+2p)!!}}\frac{h!}{(h-p)!(p-1)!} \psi_{2k+2p}(q)\bar{z}^{h-p} z^{p-1},
\end{eqnarray*}
and a shift in the summation index by one gives
$$\sum_{p=0}^{h-1} \sqrt{\frac{(2k+2p+1)!!}{(2k+2p+2)!!}}\frac{h!}{(h-p-1)!p!} \psi_{2k+2p+2}(q)\bar{z}^{h-1-p} z^{p}=h s_{e,k+1}^{h-1}$$
as required.

\hfill
$\square$

\begin{proposition}
The operators $H_t, X_t$ and $Y_t$, see \eqref{HXYholom}, preserve the span of even elements 
of the basis \eqref{hermEvenK}, \eqref{hermEven0} of the kernel 
of the symplectic Dirac operator, and the action on the basis elements 
is, for $k\in \mZ$, $k\leq -h$, given by
\begin{eqnarray}
&&H_t s^h_{e,k}=(h+2k+\frac{1}{2}) s^h_{e,k},\nonumber\\
&&X_t s^h_{e,k}=\mathrm{i}(h+k+1) s^h_{e,k+1},\nonumber\\
&&Y_t s^h_{e,k}=\mathrm{i}(k-\frac{1}{2}) s^h_{e,k-1},\quad \mbox{for } k\neq -h,\quad Y_t s^h_{e,-h}=0.
\end{eqnarray}
\end{proposition}

{\bf Proof}:
This is again a straightforward computation. For example, let us prove 
the relation $X_t s^h_{e,k}=\mathrm{i}(h+k+1) s^h_{e,k+1}$, $k\in \mN_0$:
\begin{eqnarray*}
&& \Big(\mathrm{i}\bar{z}\partial_{z}+\frac{\mathrm{i}}{4}(q-\partial_q)^2\Big)\sum_{p=0}^h \sqrt{\frac{(2k+2p-1)!!}{(2k+2p)!!}}\binom{h}{p} \psi_{2k+2p}(q)\bar{z}^{h-p} z^p\\
&& =\mathrm{i}\sum_{p=1}^{h} \sqrt{\frac{(2k+2p-1)!!}{(2k+2p)!!}}\frac{h!}{(h-p-1)!(p-1)!} \psi_{2k+2p}(q)\bar{z}^{h-p+1} z^{p-1}\\
&& +\frac{\mathrm{i}}{4}\sum_{p=0}^h \sqrt{\frac{(2k+2p-1)!!}{(2k+2p)!!}}\binom{h}{p}\times\\
&&\quad\quad\quad\quad\quad\quad\quad\times\,\, 2\sqrt{(2k+2p+1)(2k+2p+2)}\psi_{2k+2p+2}(q)\bar{z}^{h-p} z^p,
\end{eqnarray*}
and a shift by one in the summation index in the first sum gives
\begin{eqnarray*}
&& =\mathrm{i}\sum_{p=0}^h \sqrt{\frac{(2k+2p+1)!!}{(2k+2p+2)!!}}
\binom{h}{p} \psi_{2k+2p+2}(q)\bar{z}^{h-p} z^p\big(h-p+\frac{1}{2}(2k+2p+2)\big)\\
&& =\mathrm{i}(h+k+1) s^h_{e,k+1}.
\end{eqnarray*}

\hfill
$\square$

\begin{proposition}
The operators $Z_1$ and $Z_2$, see \eqref{Z1Z2zbarz}, increase the homogeneity by one 
in the variables $z,\bar{z}$ and 
preserve the even elements of the basis \eqref{hermEvenK}, \eqref{hermEven0} of the kernel 
of the symplectic Dirac operator. The operators satisfy for $k\in \mZ$, $k\leq -h$,
\begin{eqnarray}
&& Z_1 s^h_{e,k}=2(h+1)(h+k+1) s^{h+1}_{e,k},\nonumber\\
&& Z_2 s^h_{e,k}=(h+1)(2k-1) s^{h+1}_{e,k-1}.
\end{eqnarray}
\end{proposition}

{\bf Proof}:
It follows from \eqref{Z1Z2zbarz} that
\begin{eqnarray*}
Z_1&=&-\frac{1}{2} \Big( \big(q-\partial_q\big)^2 z^2+2\big( q^2 +\partial_q^2\big)z\bar{z}+\big( q+\partial_q\big)^2 \bar{z}^2\Big)\partial_z\\
&&+\bar{z}(E+1)(2E+1)+\frac{1}{2}\Big(\big( q-\partial_q \big)^2 z+\big( q+\partial_q \big)\big( q-\partial_q \big) \bar{z}\Big)(2E+1)\\
Z_2&=&-\frac{1}{2} \Big( \big(q-\partial_q\big)^2 z^2+2\big( q^2 +\partial_q^2\big)z\bar{z}+\big( q+\partial_q\big)^2 \bar{z}^2\Big)\partial_{\bar{z}}\\
&&-z(E+1)(2E+1)+\frac{1}{2}\Big(\big(q-\partial_q \big)\big( q+\partial_q \big)z+\big( q+\partial_q \big)^2 \bar{z}\Big)(2E+1)
\end{eqnarray*}
Now using \eqref{vzorceqDq}, \eqref{vzorceqDq2} and $(2E+1)s^h_{e,k}=(2h+1)s^h_{e,k}, (E+1)s^h_{e,k}=(h+1)s^h_{e,k}$,
we verify the relation $Z_2 s^h_{e,k}=(h+1)(2k-1) s^{h+1}_{e,k-1}$ for $k\in \mN_0$:
\begin{eqnarray*}
&& Z_2 \Big(\sum_{p=0}^h \sqrt{\frac{(2k+2p-1)!!}{(2k+2p)!!}}
\binom{h}{p} \psi_{2k+2p}(q)\bar{z}^{h-p} z^p\Big) 
\end{eqnarray*}
is equal to the sum of three summations:
\begin{eqnarray*}
&&= -\frac{1}{2}\sum_{p=0}^{h-1} \sqrt{\frac{(2k+2p-1)!!}{(2k+2p)!!}}\binom{h}{p} (h-p) \\
&&\quad\quad\quad\quad\quad\quad\quad 2\Big(
\sqrt{(2k+2p+1)(2k+2p+2)}\psi_{2k+2p+2}(q)\bar{z}^{h-p-1}z^{p+2}\\
&&\quad\quad\quad\quad\quad\quad\quad +(4k+4p+1)\psi_{2k+2p}(q)\bar{z}^{h-p}z^{p+1}\\
&&\quad\quad\quad\quad\quad\quad\quad +\sqrt{(2k+2p)(2k+2p-1)}\psi_{2k+2p-2}(q)\bar{z}^{h-p+1}z^{p}\Big)\\
&&-\sum_{p=0}^h \sqrt{\frac{(2k+2p-1)!!}{(2k+2p)!!}}\binom{h}{p} (h+1)(2h+1)\psi_{2k+2p}(q)\bar{z}^{h-p}z^{p+1}\\
&&+\frac{1}{2}\sum_{p=0}^h \sqrt{\frac{(2k+2p-1)!!}{(2k+2p)!!}}\binom{h}{p}(2h+1) 2\Big(
(2k+2p)\psi_{2k+2p}(q)\bar{z}^{h-p}z^{p+1}\\
&&\quad\quad\quad\quad\quad\quad\quad +\sqrt{(2k+2p)(2k+2p-1)}\psi_{2k+2p-2}(q)\bar{z}^{h-p+1}z^{p} \Big).
\end{eqnarray*}
We reorganize the sums to get the contributions to a given 
Hermite function,
\begin{eqnarray*}
&&=\sum_{p=0}^{h-1} \sqrt{\frac{(2k+2p+1)!!}{(2k+2p+2)!!}}\psi_{2k+2p+2}(q)\bar{z}^{h-p-1}z^{p+2}
\binom{h}{p}(p-h)(2k+2p+2)\\
&&+\sum_{p=0}^h \sqrt{\frac{(2k+2p-1)!!}{(2k+2p)!!}}\psi_{2k+2p}(q)\bar{z}^{h-p}z^{p+1}\binom{h}{p}\\
&&\quad\quad\quad\quad\quad\quad\quad \big((p-h)(4k+4p+1)+(2h+1)(2k+2p-h-1)\big)\\
&&+\sum_{p=0}^h \sqrt{\frac{(2k+2p-3)!!}{(2k+2p-2)!!}}\psi_{2k+2p-2}(q)\bar{z}^{h-p+1}z^{p}\binom{h}{p}(2k+2p-1)(h+p+1),
\end{eqnarray*}
and do appropriate shifts in summations and multiple 
expressions to produce the required combinatorial coefficients:
\begin{eqnarray*}
&&=\sum_{p=0}^{h+1} \sqrt{\frac{(2k+2p-3)!!}{(2k+2p-2)!!}}\psi_{2k+2p-2}(q)\bar{z}^{h-p+1}z^{p}
\frac{(h+1)!}{(h-p+1)!p!}\\
&& \frac{1}{h-1}\big(-(p-1)p(2k+2p-2)+p(p-1-h)(4k+4p-3)\\
&&+p(2h+1)(2k+2p-h-3)+(h-p+1)(2k+2p-1)(h+p+1)
\big)\\
&&= s^{h+1}_{e,k-1}(h+1)(2k-1).
\end{eqnarray*}
The remaining equalities are analogous.

\hfill
$\square$

The proof of the analogous statement for the odd part of the basis of the kernel 
of the symplectic Dirac operator is analogous to the even part
in the previous proposition and so is omitted. This result is summarized in
the next proposition.

\begin{proposition}
\begin{enumerate}
\item

The operators $\partial_z$ and $\partial_{\bar{z}}$ decrease the homogeneity in 
the variables $z,\bar{z}$
by one and preserve odd elements of the basis \eqref{hermOdd}, $k\in \mN_0$, of the kernel of the symplectic Dirac operator $D_s$:
\begin{eqnarray}
&&\partial_z s^h_{o,k} = h s^{h-1}_{o,k+1},\nonumber\\
&&\partial_{\bar{z}} s^h_{o,k} = h s^{h-1}_{o,k}.
\end{eqnarray}

\item
The operators $H_t,X_t$ and $Y_t$ preserve odd elements of the basis \eqref{hermOdd}, $k\in \mN_0$, of the kernel of the symplectic Dirac operator:
\begin{eqnarray}
&&H_t s^h_{o,k}=(h+2k+\frac{3}{2}) s^h_{o,k},\nonumber\\
&&X_t s^h_{o,k}=i(h+k+\frac{3}{2}) s^h_{o,k+1},\nonumber \\
&&Y_t s^h_{o,k}=ik s^h_{o,k-1}, \quad \mbox{for } k\neq 0, \quad Y_t s^h_{o,0}=0.
\end{eqnarray}

\item
The operators $Z_1$ and $Z_2$ increase the homogeneity by one in the 
variables $z,\bar{z}$, and map odd elements of 
the basis \eqref{hermOdd}, $k\in \mN_0$, to the elements of odd basis 
of the homogeneity plus one higher of the kernel of the symplectic 
Dirac operator,
\begin{eqnarray}
&& Z_1 s^h_{o,k}=(1+h)(2h+2k+3) s^{h+1}_{o,k},\nonumber\\
&& Z_2 s^h_{o,k}=2(h+1)k s^{h+1}_{o,k-1}.
\end{eqnarray}
\end{enumerate}
\end{proposition}

\vspace{0.5cm}
{\bf Acknowledgement:}  M. Holíková is supported by the Faculty of Mathematics and Physics (Charles University), grant SVV-2015-260227. P. Somberg is supported by the Grant Agency of Czech Republic, grant GA CR P201/12/G028.
H. De Bie is supported by the Fund for Scientific Research-Flanders (FWO-V), grant G.0116.13N.
\vspace{0.2cm}

\vspace{0.3cm}

Hendrik De Bie

Department of Mathematical Analysis, Faculty of Engineering and 
Architecture, 

Ghent University, 
Galglaan 2, 9000 Gent, Belgium

E-mail: Hendrik.DeBie@UGent.be.

\vspace{0.3cm}

Marie Hol\'ikov\'a

Department of Mathematics and Mathematical Education, Faculty of Education, 

Charles University,
Magdal\'eny Rettigov\'e 4, Praha 1, Czech Republic.

E-mail: marie.holikova@pedf.cuni.cz.

\vspace{0.3cm}

Petr Somberg

Mathematical Institute of Charles University,

Sokolovsk\'a 83, Praha 8 - Karl\'{\i}n, Czech Republic, 

E-mail: somberg@karlin.mff.cuni.cz.

\end{document}